\documentclass[10pt]{article}
\usepackage{amssymb,amsmath,textcomp}
\newtheorem{satz}{Theorem}[section]
\newtheorem{lemma}{Lemma}[section]

\newtheorem{bemerk1}{Remark}[section]
\newtheorem{bei}{Example}[section]
\newtheorem{korollar}{Corollary}[section]

\newcommand{\iR}{\mathbb{R}}
\newcommand{\iN}{\mathbb{N}}

\newcommand{\iC}{\mathbb{C}}

\newcommand{\oH}{\hspace*{0.39em}\raisebox{0.6ex}{\textdegree}\hspace{-0.72em}H}
\DeclareMathOperator*{\esup}{ess\,sup}
\DeclareMathOperator*{\einf}{ess\,inf}

\begin{document}
\begin{center}
{\bf\Large Optimal decay estimates for time-fractional and other non-local subdiffusion equations
via energy methods}
\end{center}
\vspace{0.7em}
\begin{center}
Vicente Vergara\footnote{The first author was partially supported by FONDECYT grant 1110033.} and Rico Zacher\footnote{The second author was supported by a Heisenberg fellowship of the German Research Foundation (DFG).}
\end{center}
\begin{abstract}
We prove sharp estimates for the decay in time of solutions to a rather general class of non-local in time subdiffusion
equations on a bounded domain subject to a homogeneous Dirichlet boundary condition. Important special cases are
the time-fractional and ultraslow diffusion equation, which have seen much interest during the last years, mostly due to their applications in the modeling of anomalous diffusion. We study the case where the equation is in divergence form with bounded measurable coefficients. Our proofs rely on energy estimates and make use of a new and powerful inequality for integro-differential operators of the form $\partial_t (k\ast \cdot)$. The results can be generalized to certain quasilinear equations. We illustrate this by looking at the time-fractional $p$-Laplace and porous medium equation. Here it turns out
that the decay behaviour is markedly different from that in the classical parabolic case.

\end{abstract}
\vspace{0.7em}
\begin{center}
\end{center}

\section{Introduction and main results}
Let $\Omega\subset \iR^N$ be a bounded domain. We are interested in
the long-time behaviour of solutions to the non-local in time diffusion equation
\begin{equation} \label{lindiff}
\partial_t \big(k\ast [u-u_0]\big)-\mbox{div}\,\big(A(t,x)Du\big)=0,\quad
t>0,\,x\in \Omega,
\end{equation}
subject to the Dirichlet boundary condition
\begin{equation} \label{lindiff2}
u|_{\partial \Omega}=0,\quad t>0,\,x\in \partial\Omega.
\end{equation}
Here $u_0=u_0(x)$ plays the role of the initial datum for $u$, that is
\begin{equation} \label{lindiff3}
u|_{t=0}=u_0,\,x\in \partial\Omega.
\end{equation}
The kernel $k\in L_{1,\,loc}(\iR_+)$ is given, and $k\ast v$ denotes the
convolution on the positive halfline $\iR_+:=[0,\infty)$ w.r.t.\ the time variable,
that is
$(k\ast v)(t)=\int_0^t k(t-\tau)v(\tau)\,d\tau$, $t\ge 0$.
We assume that $k$ is a kernel of type $\mathcal{PC}$, by which we mean that the following condition is satisfied.
\begin{itemize}
\item [{\bf ($\mathcal{PC}$)}] $k\in L_{1,\,loc}(\iR_+)$ is nonnegative and nonincreasing, and there exists a kernel $l\in L_{1,\,loc}(\iR_+)$ such that
$k\ast l=1$ on $(0,\infty)$.
\end{itemize}
We also write $(k,l)\in \mathcal{PC}$ in this situation. From $(k,l)\in {\cal PC}$ it follows that $l$ is completely
positive, see e.g. Theorem 2.2 in \cite{CN}, in particular $l$ is nonnegative.

Concerning the
coefficients $A=(a_{ij})$ we merely assume measurability,
boundedness, and a uniform parabolicity condition, that is we assume
that
\begin{itemize}
\item [{\bf ($\mathcal{H}$)}] $A\in L_\infty((0,T)\times \Omega;\iR^{N\times
N})$ for all $T>0$, and $\exists \nu>0$ such that
\[
\big(A(t,x)\xi|\xi\big)\ge \nu|\xi|^2,\quad \mbox{for
a.a.}\;(t,x)\in (0,\infty)\times \Omega,\,\mbox{and all}\,\xi\in
\iR^N.
\]
\end{itemize}

An important example for a pair $(k,l)\in \mathcal{PC}$ is given by $(k,l)=(g_{1-\alpha},g_\alpha)$ with $\alpha\in(0,1)$, where
$g_\beta$ denotes the standard kernel
\[
g_\beta(t)=\,\frac{t^{\beta-1}}{\Gamma(\beta)}\,,\quad
t>0,\quad\beta>0.
\]
In this case, (\ref{lindiff}) is an equation of {\em fractional time order}
$\alpha\in (0,1)$, often called {\em time-fractional diffusion equation} for $A=\nu I$; here the term $\partial_t(k\ast v)$ becomes the classical Riemann-Liouville fractional derivative
$\partial_t^\alpha v$ of the (sufficiently smooth) function $v$, see e.g.\ \cite{KST}.

Another interesting example is given by the pair
\begin{equation} \label{ultrapair}
k(t)=\int_0^1 g_\beta(t)\,d\beta,\quad l(t)=\int_0^\infty \,\frac{e^{-st}}{1+s}\,{ds},\quad t>0.
\end{equation}
In this case the operator $\partial_t(k\ast \cdot)$ is a so-called operator of
{\em distributed order}, see e.g.\ \cite{Koch08}. More examples will be discussed in Section \ref{kernelexamples} below.

Equations of the form (\ref{lindiff}) appear in mathematical physics in the context
of anomalous diffusion processes, see e.g.\ \cite{Koch08}, \cite{Koch11}, \cite{Metz}, \cite{Uch}. Let us consider for a moment the situation where $\Omega=\iR^N$ and $A=I$. Denote by $Z(t,x)$ the fundamental solution of (\ref{lindiff}) with $Z(0,x)=\delta(x)$. For $k$ as in the previous examples it is
known that $Z(t,\cdot)$ is a probability density function for all $t>0$. An important quantity that describes how fast particles diffuse and which can be measured in experiments is the {\em mean
square displacement} which is defined as
\[
m(t)=\int_{\iR^N}|x|^2 Z(t,x)\,dx,\quad t>0.
\]
In the case of the classical diffusion equation (i.e.\ $\alpha=1$) $m(t)=ct$, $t>0$
with some constant $c>0$. In the time-fractional diffusion case (i.e.\ the first
example) one observes that $m(t)=ct^\alpha$ (cf.\ \cite{Metz}), which shows that the diffusion is slower than in the classical case of Brownian motion. In the second example,
$m(t)$ behaves like $c\log t$ for $t\to \infty$, see \cite{Koch08}. In this case
(\ref{lindiff}) describes a so-called {\em ultraslow diffusion} process.

Another context where equations of the form (\ref{lindiff}) and nonlinear variants of them arise is the modelling
of dynamic processes in materials with {\em memory}. Examples are given by the theory of heat
conduction with memory, see \cite{JanI} and the references therein, and the
diffusion of fluids in porous media with memory, see \cite{CapuFlow}, \cite{JakuDiss}.

One of the main goals of this paper is to prove sharp decay estimates for suitably defined solutions
to (\ref{lindiff})--(\ref{lindiff3}). Among others, it turns out that the $L_2(\Omega)$-norm of $u$ decays algebraically like $ct^{-\alpha}$ in the first example, whereas in the ultraslow diffusion example it behaves like $c(\log t)^{-1}$ for $t\to \infty$. Recall that in the classical case one has an exponential decay
(as $\Omega$ is assumed to be bounded). These decay rates reflect, like the mean square displacement, the different degrees of slowness
of diffusion in these examples.

To motivate our first main result, let us consider the special case $A(t,x)=I$, that is the equation
\begin{equation} \label{deltaequ}
\partial_t \big(k\ast [u-u_0]\big)-\Delta u=0,\quad
t>0,\,x\in \Omega,
\end{equation}
together with (\ref{lindiff2}) and (\ref{lindiff3}). Assume that $u_0\in L_2(\Omega)$. Let $\{\phi_n\}_{n=1}^\infty
\subset \oH^1_2(\Omega)$ be an orthonormal basis of $L_2(\Omega)$ consisting of eigenfunctions of the negative
Dirichlet Laplacian with eigenvalues $\lambda_n>0$,
$n\in \iN$. Denote by $\lambda_1$ the smallest such eigenvalue. For $\mu\ge 0$ define the relaxation function $s_\mu$ on $[0,\infty)$ as the solution of the Volterra
equation
\begin{equation} \label{smudef}
s_\mu(t)+\mu(l\ast s_\mu)(t)=1,\quad t\ge 0.
\end{equation}
Note that $s_0\equiv 1$ and that (\ref{smudef}) is equivalent to the integro-differential equation
\[
\frac{d}{dt}\,\left(k\ast [s_\mu-1]\right)(t)+\mu s_\mu(t)=0,\quad t>0,\quad s_\mu(0)=1.
\]
It is known that the assumption $(k,l)\in \mathcal{PC}$ implies that $s_\mu$ is nonnegative, nonincreasing, and that $s_\mu\in H^1_{1,\,loc}(\iR_+)$; moreover
$\partial_\mu s_\mu(t)\le 0$, see e.g.\ Pr\"uss \cite{JanI}.
The solution $u$ can now
be represented via Fourier series as
\begin{equation} \label{uformel}
u(t,x)=\sum_{n=1}^\infty s_{\lambda_n}(t)\,(u_0|\phi_n)\phi_n(x),\quad t\ge 0,\,x\in \Omega,
\end{equation}
where $(\cdot|\cdot)$ denotes the standard inner product in $L_2(\Omega)$, cf.\ \cite[Theorem 4.1]{NSY} for the special case $k=g_{1-\alpha}$. By Parseval's identity and since $\partial_\mu s_\mu\le 0$, it follows from (\ref{uformel}) that
\begin{align*}
|u(t,\cdot)|_{L_2(\Omega)}^2 & =\sum_{n=1}^\infty s_{\lambda_n}^2(t)\,|(u_0|\phi_n)|^2\\
& \le s_{\lambda_1}^2(t) \sum_{n=1}^\infty|(u_0|\phi_n)|^2\\
& = s_{\lambda_1}^2(t) |u_0|_{L_2(\Omega)}^2,
\end{align*}
and thus
\begin{equation} \label{udecay}
|u(t,\cdot)|_{L_2(\Omega)}\le s_{\lambda_1}(t)|u_0|_{L_2(\Omega)},\quad t\ge 0.
\end{equation}
This decay estimate is optimal as the example $u_0=\phi_1$ with solution $u(t,x)
=s_{\lambda_1}(t)\phi_1(x)$ shows. To our knowledge the estimate (\ref{udecay})
for solutions of (\ref{deltaequ}) seems to be new in the case of a general kernel $k$ enjoying property $(\mathcal{PC})$. The special case $k=g_{1-\alpha}$ can already be found in \cite{MNV} and \cite[Corollary 4.1]{NSY}.
Concerning the long-time behaviour of solutions to abstract linear and nonlinear Volterra equations we also refer to \cite{AP}, \cite{CN}, \cite{JanI}, and \cite{SaVe13} and the references given therein.

One of the purposes of this paper is to generalize the decay estimate (\ref{udecay}) to the weak setting
with an operator in divergence form as described as above.

Let $u_0\in L_2(\Omega)$.
We say that a function $u$ is a {\em weak solution (subsolution,
supersolution)} of
(\ref{lindiff})--(\ref{lindiff3}) on $\Omega_T:=(0,T)\times \Omega$ if $u$
belongs to the space
\begin{align*}
V(T):=\{&\,v\in
L_2([0,T];\oH^1_2(\Omega))\;
\mbox{such that}\;\\
&\;\;k\ast v\in C([0,T];L_2(\Omega)),
\;\mbox{and}\;(k\ast v)|_{t=0}=0\},
\end{align*}
and for any nonnegative test function
\[
\eta\in \oH^{1,1}_2(\Omega_T)=H^1_2([0,T];L_2(\Omega))\cap
L_2([0,T];\oH^1_2(\Omega))\quad\quad
\Big(\oH^1_2(\Omega):=\overline{C_0^\infty(\Omega)}\,{}^{H^1_2(\Omega)}\Big)
\]
with $\eta|_{t=T}=0$ there holds
\[
\int_{0}^{T} \int_\Omega \Big(-\eta_t \big(k\ast [u-u_0]\big)+
(ADu|D \eta)\Big)\,dx\,dt=\,(\le,\,\ge)\,0.
\]
We say that a
function $u:(0,\infty)\times \Omega \rightarrow \iR$ is a {\em global weak solution (subsolution, supersolution)} of
(\ref{lindiff})--(\ref{lindiff3}) if for any $T>0$ the restriction
$u|_{(0,T)\times \Omega}$ is a weak solution (subsolution, supersolution) of
(\ref{lindiff})--(\ref{lindiff3}) on $(0,T)\times \Omega$.

We remark that existence and uniqueness of a global weak solution to
(\ref{lindiff})--(\ref{lindiff3}) under the above assumptions
follow from the results in \cite{ZWH}. Observe that $u\in V(T)$ does not imply
$u\in C([0,T];L_2(\Omega))$ in general, so it is not so clear how to interpret the
initial condition. However, once one knows that the functions $u$ and $k\ast
(u-u_0)$ are sufficiently smooth, then $u|_{t=0}=u_0$ is satisfied in an appropriate sense (see \cite{ZWH}). We further mention that for any weak solution of (\ref{lindiff})--(\ref{lindiff3}) on $(0,T)\times \Omega$ we also have
$\frac{d}{dt}(k\ast (u-u_0))\in L_2([0,T];H^{-1}_2(\Omega))$, where the time derivative has to be understood in the generalized sense and $H^{-1}_2(\Omega))$
denotes the dual space of $\oH^1_2(\Omega)$, see \cite{ZWH}.

Notice also that under the above assumptions the weak maximum principle is valid and takes the same form as in the classical parabolic case (see \cite{Za}). Thus the global weak solution $u$ of (\ref{lindiff})--(\ref{lindiff3}) satisfies
\begin{equation} \label{maxprin}
\einf_{\Omega}u_0\le u(t,x)\le
\esup_{\Omega}u_0,\quad \mbox{for
a.a.}\;(t,x)\in (0,\infty)\times \Omega,
\end{equation}
provided $u_0\in L_\infty(\Omega)$. We also refer to \cite{Luch} for a different proof of the maximum principle in the more special situation of strong solutions to the time-fractional diffusion equation. We further remark that in the special case $k=g_{1-\alpha}$ with $\alpha\in (0,1)$ H\"older continuity of the weak solution to (\ref{lindiff})--(\ref{lindiff3}) with $u_0\in L_\infty(\Omega)$ has
been established recently in \cite{Za1}, see also \cite{ZaGS}.

Denoting by $y_+$ and $y_-:=[-y]_+$ the positive and negative part, respectively, of $y\in \iR$, our first main result concerning (\ref{lindiff})--(\ref{lindiff3}) reads as follows.
\begin{satz} \label{result1}
Let $u_0\in L_2(\Omega)$ and suppose that the conditions (H) and  ($\mathcal{PC}$) are satisfied.
Then for any global weak subsolution (supersolution) $u$ of (\ref{lindiff})--(\ref{lindiff3}), there holds
\[
\big|u_{+\,(-)}(t,\cdot)\big|_{L_2(\Omega)}\le s_{\nu \lambda_1}(t)\,\big|[u_0]_{+\,(-)}\big|_{L_2(\Omega)},\quad
\mbox{a.a.}\,t>0.
\]
\end{satz}
As a direct consequence we obtain
\begin{korollar} \label{korresult1}
Let $u_0\in L_2(\Omega)$ and assume that the conditions (H) and  ($\mathcal{PC}$) are fulfilled.
Then the global weak solution $u$ of (\ref{lindiff})--(\ref{lindiff3}) satisfies the estimate
\begin{equation} \label{korresest}
|u(t,\cdot)|_{L_2(\Omega)}\le s_{\nu \lambda_1}(t)\,|u_0|_{L_2(\Omega)},\quad
\mbox{a.a.}\,t>0.
\end{equation}
\end{korollar}
These decay estimates are again optimal as the special case $A=\nu I$ shows, in fact specializing further to $\nu=1$ we recover the estimate (\ref{udecay}).

We would like to point out that even though (\ref{lindiff}) is linear, Theorem \ref{result1} and Corollary \ref{korresult1} are {\em nonlinear} results. For example one could think of $A(t,x)=A_0(t,x,u(t,x))$ with some appropriate
nonlinear function $A_0$. It is further possible to extend these results without
much effort to quasilinear equations of the form
\[
\partial_t \big(k\ast [u-u_0]\big)-\mbox{div}\,A(t,x,u,Du)=0,\quad
t>0,\,x\in \Omega,
\]
where $A$ satisfies suitable measurability and structure conditions, in particular
\[(A(t,x,u,Du)|Du)\ge \nu|Du|^2\quad\mbox{with some}\; \nu>0.\]
To illustrate this aspect, we mention a quasilinear time-fractional problem with Dirichlet boundary condition which has been studied recently in \cite{ZaGS}, where $A(t,x)=A(u(t,x))$. There it was shown that the $L_2$-norm of the solution decays
at least like $t^{-\alpha/2}$. Applying Corollary \ref{korresult1} not only improves this estimate, but also provides the optimal result
which says that $|u(t,\cdot)|_2$ decays like $t^{-\alpha}$.

The proof of Theorem \ref{result1} is based on suitable energy estimates and a new and extremely useful
inequality for integro-differential operators of the form $\partial_t(k\ast \cdot)$ we will refer to as the {\em $L_p$-norm inequality for $\partial_t(k\ast \cdot)$}, see Section \ref{SecLpIn} below. For $p=2$ and $k=g_{1-\alpha}$ with $\alpha\in (0,1)$ it takes the form
\begin{equation} \label{introungl}
|u(t)|_{L_2(\Omega)}\partial_t^\alpha\Big( |u(\cdot)|_{L_2(\Omega)}-|u_0|_{L_2(\Omega)}\Big)(t)\le \int_\Omega u(t,x) \,\partial_t^\alpha(u-u_0)(t,x)\,dx,\quad \mbox{a.a.}\;t\in (0,T),
\end{equation}
for all $u_0\in L_2(\Omega)$ and all sufficiently smooth functions $u:[0,T]\times \Omega\rightarrow \iR$. Once (\ref{introungl}) is known, it is quite straightforward to prove the desired decay rate for {\em solutions} in the time-fractional case by formal estimates. In fact,
testing the PDE with $u$, integrating over $\Omega$, and using $A\ge \nu I$ as well as Poincar\'e's inequality gives
\[
\int_\Omega u \,\partial_t^\alpha(u-u_0)\,dx +\nu \lambda_1 \int_\Omega |u|^2\,dx\le 0,\quad t>0.
\]
By (\ref{introungl}) this implies
\[
|u(t)|_{L_2(\Omega)}\partial_t^\alpha\Big( |u(\cdot)|_{L_2(\Omega)}-|u_0|_{L_2(\Omega)}\Big)(t)+\nu\lambda_1|u(t)|_{L_2(\Omega)}^2\le 0,\quad t>0.
\]
Assuming $|u(t)|_{L_2(\Omega)}>0$ we thus obtain the fractional differential inequality
\[
\partial_t^\alpha\big(|u|_{L_2(\Omega)}-|u_0|_{L_2(\Omega)}\big)(t)+\nu\lambda_1|u(t)|_{L_2(\Omega)}\le 0,\quad t>0,
\]
which implies (\ref{korresest}), by a comparison principle argument (see Section \ref{VoltHilf} below). To give a rigorous proof of Theorem \ref{result1}, which is on {\em sub-} and {\em supersolutions} in the {\em weak} setting, requires much more effort. In particular the problem has to be regularized in time suitably in order to justify the application of the so-called fundamental identity for operators of the form $\partial_t(k\ast \cdot)$ (see (\ref{fundidentity}) below), which is the basic
tool for deriving a priori estimates for equations of the form (\ref{lindiff}) (cf.\ \cite{Za}) and also the key ingredient in the proof of the $L_p$-norm inequality.

Our techniques also apply to other types of non-local in time subdiffusion equations. In the present paper we also consider the time-fractional $p$-Laplace equation
\begin{align*}
\partial_t^\alpha(u-u_0)-\Delta_p u & =0\quad \mbox{in}\;\iR_+\times \Omega,
\nonumber\\
u|_{\partial \Omega} & = 0 \quad \mbox{at}\; \iR_+\times \partial \Omega,
\\
u|_{t=0} & = u_0 \quad \mbox{in}\; \Omega,\nonumber
\end{align*}
where $\alpha\in (0,1)$ and $1<p<\infty$. It is well known that in the classical case $\alpha=1$, solutions decay algebraically as $t\to \infty$ if $p>2$,
whereas for $p<2$ one has the phenomenon of extinction in finite time (\cite{DBUV}). It turns out that in the time-fractional case solutions decay algebraically like $t^{-\frac{\alpha}{p-1}}$ in the whole range of $p$, see Theorem \ref{plaplaceresult}; our results indicate that extinction in finite time does not occur anymore, at least if $p>\frac{2N}{N+2}$. This interesting phenomenon is due to the slowness of the diffusion in the case $\alpha<1$. A corresponding result can be shown for the time-fractional porous medium equation,
see Theorem \ref{mlaplaceresult}.

The paper is organized as follows. In Section 2 we collect some preliminary results on operators of the type $\partial_t(k\ast \cdot)$ and related Volterra integral equations. In particular we prove rather general comparison results for these equations, which seem
to be new and are interesting in its own right. Section 3 is devoted to the $L_p$-norm inequality and several variants of it. In Section 4 we derive a subsolution inequality for the positive part of a subsolution to (\ref{lindiff})--(\ref{lindiff3}). This is an important step in our proof of Theorem \ref{result1}, which is completed in Section 5. In Section 6 we illustrate our results by looking at several
examples of pairs of kernels $(k,l)\in \mathcal{PC}$. We will see that this class of kernels is quite rich in that the solutions may exhibit a very different kind of decay like e.g.\ exponential, algebraic, or logarithmic decay. We obtain sharp decay rates, making use of results on the Laplace transform and the Karamata-Feller Tauberian theorem.
In Section 7 we discuss the asymptotic behaviour of the solution to the nonlinear time-fractional differential equation
 \[
\partial_t^\alpha(u-u_0)+\nu u^\gamma=0,\;\;t\ge 0,\quad\; u(0)=u_0>0,
\]
 where $\alpha\in (0,1)$ and $\nu,\gamma>0$. This equation is fundamental for deriving optimal decay estimates
 for the quasilinear problems studied in the last two sections, Sections 8 and 9, which are concerned with the time-fractional $p$-Laplace and porous medium equation, respectively.
\section{Preliminaries}
\subsection{Regularization of the kernel} Let $(k,l)\in \mathcal{PC}$. For $\mu>0$ let $h_\mu\in L_{1,loc}(\iR_+)$ denote the resolvent kernel associated
with $\mu l$, that is we have
\begin{equation} \label{hndef}
h_\mu(t)+\mu(h_\mu\ast l)(t)=\mu l(t),\quad t>0,\,\mu>0.
\end{equation}
Note that $h_\mu=-\dot{s}_\mu \in L_{1,\,loc}(\iR_+)$, in particular $h_\mu$ is nonnegative.
It is well-known that for any $f\in L_p([0,T])$, $1\le p<\infty$, there holds
$h_n\ast f\rightarrow f$ in $L_p([0,T])$ as $n\rightarrow \infty$, see e.g.\ \cite{Za}.

For $\mu>0$ we set
\begin{equation} \label{kndef}
k_{\mu}=k\ast h_\mu.
\end{equation}
It is known (see e.g.\ \cite{Za}) that $k_\mu=\mu s_\mu$, $\mu>0$, and thus
the kernels $k_\mu$ are also nonnegative
and nonincreasing, and they belong to $H_1^1([0,T])$ for any $T>0$.

The following lemma, which can be found in \cite{Za}, provides an equivalent weak formulation where
the singular kernel $k$ appearing in the integro-differential operator w.r.t.\ time is approximated by a more regular kernel.
\begin{lemma} \label{regWF}
Let $u_0\in L_2(\Omega)$ and suppose that the conditions (H) and  ($\mathcal{PC}$) are satisfied. Then
$u\in V(T)$ is a weak solution (subsolution, supersolution) of (\ref{lindiff})--(\ref{lindiff3}) on $\Omega_T$
if and only if for any nonnegative
function $\psi\in \oH^1_2(\Omega)$ there holds
\begin{align} \label{regcond}
\int_\Omega
\Big(\psi\partial_t [k_n\ast(u-u_0)]+\big(h_n\ast[ADu]|D\psi\big)\Big)\,dx=\,(\le,\,\ge)\,0,\quad
\mbox{a.a.}\,t\in (0,T),\,n\in \iN.
\end{align}
\end{lemma}
\subsection{The fundamental identity} We next state a fundamental identity for integro-differential
operators of the form $\partial_t(k\ast \cdot)$, cf.\ also \cite{Za},
\cite{Za1}. It can be viewed as the analogue to the chain rule
$(H(u))'=H'(u)u'$.
\begin{lemma} \label{FILemma}
Let $T>0$ and $U$ be an open subset of $\iR$. Let further $k\in
H^1_1([0,T])$, $H\in C^1(U)$, and $u\in L_1([0,T])$ with $u(t)\in U$
for a.a. $t\in (0,T)$. Suppose that the functions $H(u)$, $H'(u)u$,
and $H'(u)(\dot{k}\ast u)$ belong to $L_1([0,T])$ (which is the case
if, e.g., $u\in L_\infty([0,T])$). Then we have for a.a. $t\in
(0,T)$,
\begin{align} \label{fundidentity}
H'(u(t))&\frac{d}{dt}\,(k \ast u)(t) =\;\frac{d}{dt}\,\big(k\ast
H(u)\big)(t)+
\Big(-H(u(t))+H'(u(t))u(t)\Big)k(t) \nonumber\\
 & +\int_0^t
\Big(H(u(t-s))-H(u(t))-H'(u(t))[u(t-s)-u(t)]\Big)[-\dot{k}(s)]\,ds.
\end{align}
\end{lemma}
The lemma follows from a straightforward computation. In particular
identity (\ref{fundidentity}) applies to the regularized operator $u\mapsto \partial_t( k_n\ast u)$ from above. We remark that an integrated
version of (\ref{fundidentity}) can be found in \cite[Lemma
18.4.1]{GLS}. Observe that (\ref{fundidentity}) remains valid for singular kernels $k$, like e.g.\ $k=g_{1-\alpha}$ with
$\alpha\in(0,1)$, provided that $u$ is sufficiently smooth.

The special case $H(y)=\frac{1}{2}y^2$ extends to the Hilbert space setting. The following lemma can be found
in \cite{VZ}.
\begin{lemma} \label{fundlemma1}
Let $T>0$ and ${\cal H}$ be a real Hilbert space with inner product $\langle\cdot,\cdot\rangle_{{\cal H}}$. Then for any
$k\in H^1_1([0,T])$ and any $u\in L_2([0,T];{\cal H})$ there holds
\begin{align*}
\Big\langle u(t),\frac{d}{dt}\,(k\ast u)(t)\Big\rangle_{{\cal H}} =&\;
\frac{1}{2}\,\frac{d}{dt}\,(k\ast
|u|_{\cal H}^2)(t)+\frac{1}{2}\,k(t)|u(t)|_{\cal H}^2 \nonumber\\
&\;+\,\frac{1}{2}\,\int_0^t |u(t)-u(t-s)|_{\cal H}^2[-\dot{k}(s)]\,ds,\quad
\mbox{a.a.}\;t\in(0,T).
\end{align*}
\end{lemma}

${}$

\subsection{Auxiliary results on Volterra equations} \label{VoltHilf}
\begin{lemma} \label{contdep}
Let $(k,l)\in \mathcal{PC}$, $\mu> 0$, and $T>0$. Let $f_n\in L_1([0,T])$, $n\in \iN$,
and denote by $v_n\in L_1([0,T])$ the solution of
\[ \partial_t(k_n\ast w)(t)+\mu w(t)=f_n(t),\quad \mbox{a.a.}\;t\in (0,T),
\]
where $k_n$ is defined as in (\ref{kndef}). Suppose that $f_n\rightarrow f$ in
$L_1([0,T])$. Then $v_n\rightarrow v$ in $L_1([0,T])$, where $v$ solves the equation
\[
v(t)+\mu(l\ast v)(t)=(l\ast f)(t),\quad \mbox{a.a.}\;t\in (0,T).
\]
\end{lemma}
{\em Proof.} a) Note first that for any $w,\tilde{f}\in L_1([0,T])$,
\[
\partial_t(k_n\ast w)=\tilde{f},\quad \mbox{a.e. in}\;(0,T)
\]
is equivalent to
\[
w=\frac{1}{n}\,\tilde{f}+\,l\ast \tilde{f},\quad \mbox{a.e. in}\;(0,T).
\]
This can be easily seen, e.g., with the aid of the Laplace transform.

b) Using a) with $\tilde{f}=f_n-\mu v_n$ we see that
\[
v_n+\mu l\ast v_n=\frac{1}{n}f_n+l\ast f_n-\frac{\mu}{n}v_n=:F_n,
\]
which gives
\begin{equation} \label{vngl}
v_n=F_n-h_\mu\ast F_n.
\end{equation}
On the other hand we have
\begin{equation} \label{vgl}
v=l\ast f-h_\mu \ast l\ast f.
\end{equation}

c) Subtracting (\ref{vgl}) from (\ref{vngl}) and taking the $L_1([0,T])$-norm we obtain
\begin{align*}
|v_n-v|_1\le &\,|l\ast(f_n-f)|_1+\frac{1}{n}\,\big(|f_n|_1+|h_\mu\ast f_n|_1\big)
+|h_\mu\ast l\ast(f_n-f)|_1\\
&\,+\frac{\mu}{n}|v_n-v|_1+\frac{\mu}{n}|v|_1+\frac{\mu}{n}|h_\mu\ast(v_n-v)|_1+
\frac{\mu}{n}|h_\mu\ast v|_1.
\end{align*}
We may now use Young's inequality for convolutions and absorb, for large $n$, the terms on the right-hand side that involve
$(v_n-v)$-factors to conclude that $v_n\rightarrow v$ in $L_1([0,T])$. \hfill $\square$
\begin{lemma} \label{compkngl}
Let $(k,l)\in \mathcal{PC}$, $\mu\ge 0$, and $T>0$. Let $n\in \iN$ and $k_n$ be defined as in (\ref{kndef}). Suppose further that $u,v,f\in L_1([0,T])$ are such that
\begin{align*}
\partial_t(k_n\ast u)+\mu u\le & \,f,\quad \mbox{a.e. in}\; (0,T);\\
\partial_t(k_n\ast v)+\mu v\ge & \,f,\quad \mbox{a.e. in}\; (0,T).
\end{align*}
Then $u\le v$ a.e.\ in $(0,T)$.
\end{lemma}
{\em Proof.} From the assumptions it follows immediately that
\[
\partial_t(k_n\ast [u-v])+\mu (u-v)\le 0,\quad \mbox{a.e. in}\; (0,T).
\]
We multiply this inequality with $(u-v)_+$ and apply the fundamental identity
(\ref{fundidentity}) with $H(y)=\frac{1}{2}(y_+)^2$. This gives
\begin{equation} \label{compu1}
\frac{1}{2}\,\partial_t(k_n\ast [u-v]_+^2)+\mu [u-v]_+^2\le 0,\quad \mbox{a.e. in}\; (0,T).
\end{equation}
Next, we apply the positive operator $(\frac{1}{n}+l\ast)$ to (\ref{compu1}),
thereby getting (cf.\ step a) in the proof of Lemma \ref{contdep}) that
\[
(u-v)_+^2+\mu(\frac{1}{n}+l\ast)(u-v)_+^2\le 0,\quad \mbox{a.e. in}\; (0,T),
\]
which in turn implies the assertion. \hfill $\square$

${}$

\noindent We also need the following nonlinear comparison result with a singular kernel.
\begin{lemma} \label{compstrong}
Let $T>0$, $(k,l)\in \mathcal{PC}$, and $f\in C(\iR)$. Assume that $f$ is nondecreasing. Suppose that $v,\,w\in H^1_1([0,T])$ satisfy $v(0)\le w(0)$ and
\begin{align*}
\partial_t\big(k\ast [v-v(0)]\big)+f(v) &\le 0,\quad \mbox{a.a.}\; t\in (0,T),\\
\partial_t\big(k\ast [w-w(0)]\big)+f(w) &\ge 0,\quad \mbox{a.a.}\; t\in (0,T).
\end{align*}
Then $v(t)\le w(t)$ for all $t\in [0,T]$.
\end{lemma}
{\em Proof.} Subtracting the second from the first inequality yields
\[
\partial_t\big(k\ast [v-w]\big)+f(v)-f(w)\le \big(v(0)-w(0)\big)k(t)\le 0.
\]
We convolve this inequality with the positive kernel $h_n$ with $n\in \iN$ (see (\ref{hndef}) for its definition). Since $k\ast(v-w)(0)=0$ and using $k_n=h_n\ast k$ we obtain
\[
 \partial_t\big(k_n\ast [v-w]\big)+h_n\ast \big(f(v)-f(w)\big)\le 0.
\]
Next, we multiply by $(v-w)_+$ and apply the fundamental identity as in the previous proof to get
\[
\frac{1}{2}\,\partial_t\big(k_n\ast (v-w)_+^2\big)+\big[h_n\ast \big(f(v)-f(w)\big)\big] (v-w)_+\le 0.
\]
Convolving with the positive kernel $l$ then yields
\[
\frac{1}{2}\,h_n\ast (v-w)_+^2+l\ast \left(\big[h_n\ast \big(f(v)-f(w)\big)\big]\, (v-w)_+\right)\le 0.
\]
We now send $n\rightarrow \infty$ and select an appropriate subsequence, if necessary, to infer that
\begin{equation} \label{laststep}
\frac{1}{2}\,(v-w)_+^2+l\ast \left(\big[f(v)-f(w)\big]\, (v-w)_+\right)\le 0.
\end{equation}
Since $f$ is nondecreasing, the second term is nonnegative, and thus (\ref{laststep})
implies $(v-w)_+=0$, that is $v\le w$ in $[0,T]$. \hfill $\square$
\begin{bemerk1} \label{vergleichweak}
The last result remains true for {\em weak} sub- and supersolutions, that is
for $v,w\in L_1([0,T])$ with $k\ast v \in C([0,T])$, $k\ast w\in C([0,T])$, $f(v),f(w)\in L_1([0,T])$, and for $v_0\le w_0$ satisfying
\begin{align*}
\int_0^T \left(-\dot{\varphi}\big(k\ast [v-v_0]\big)+\varphi f(v)\right)\,dt &\le 0,\\
\int_0^T \left(-\dot{\varphi}\big(k\ast [w-w_0]\big)+\varphi f(w)\right)\,dt &\ge 0,
\end{align*}
for all nonnegative $\varphi\in H^1_1([0,T])$ with $\varphi(T)=0$. In this situation one has $v(t)\le w(t)$ for a.a.\ $t\in (0,T)$.
\end{bemerk1}
\subsection{A version of the Karamata-Feller Tauberian theorem}
The asymptotic behaviour of a function $w(t)$ as $t\to \infty$ can be determined, under suitable conditions, by looking
at the behaviour of its Laplace transform $\hat{w}(z)$ as $z\to 0$, and vice versa. An important situation where such a correspondence holds is described by the Karameter-Feller Tauberian theorem. We state a special case of it, which suffices for our purposes. See the monograph \cite{Feller} for a more general version and proofs.
\begin{satz} \label{Karamata}
Let $L:(0,\infty)\to (0,\infty)$ be a function that is {\em slowly varying at $\infty$}, that is, for every fixed
$x>0$ we have $L(tx)/L(t)\to 1$ as $t\to \infty$. Let $\beta>0$ and $w:(0,\infty)\rightarrow \iR$ be a monotone function whose Laplace transform $\hat{w}(z)$ exists for all $z\in \iC_+:=\{\lambda\in \iC:\, \mbox{Re}\,\lambda>0\}$. Then
\[
\hat{w}(z) \sim \,\frac{1}{z^\beta}\,L\left(\frac{1}{z}\right)\quad \mbox{as}\;z\to 0\quad\; \mbox{if and only if} \quad \;
w(t)\sim g_\beta(t) L(t)\quad \mbox{as}\;t\to \infty.
\]
Here the approaches are on the positive real axis and the notation $f(t)\sim g(t)$ as $t\to t_*$ means that
$\lim_{t\to t_*} f(t)/g(t)=1$.
\end{satz}
\section{The $L_p$-norm inequality} \label{SecLpIn}
The following inequality seems to be new and is the key inequality to obtain, among others, sharp $L_p$-norm decay estimates for various types of linear and nonlinear integro-differential equations
involving an operator $\partial_t(k\ast \cdot)$.

In what follows we use the convention that $y|y|^{-\beta}$ is equal to $0$ for
$y=0$ whenever $\beta\in (0,1)$.
\begin{lemma} \label{lpnorminequ}
Let $1<p<\infty$, $T>0$, and $\Omega$ be an arbitrary measurable subset of $\iR^N$. Let further $k\in
H^1_1([0,T])$ be nonnegative and nonincreasing. Then for any
$u\in L_p([0,T];L_p(\Omega))$ there holds
\[
|u(t)|_{L_p(\Omega)}^{p-1}\partial_t\big(k\ast |u(\cdot)|_{L_p(\Omega)}\big)(t)\le \int_\Omega |u|^{p-2}u \,\partial_t(k\ast u)(t)\,dx,\quad \mbox{a.a.}\;t\in (0,T).
\]
\end{lemma}
{\em Proof.} Let $H(y)=\frac{1}{p}\,|y|^p,\,y\in \iR$. Then $H'(y)=|y|^{p-2}y$ and thus by the fundamental identity, Fubini's theorem, and H\"older's inequality we have
\begin{align*}
\int_\Omega |u|^{p-2}u \partial_t(k\ast u)\,dx & =  \,\frac{1}{p}\,\int_\Omega \partial_t(k\ast|u|^p)\,dx+\int_\Omega \big(-\frac{1}{p}\,|u|^p+|u|^p\big)k(t)\,dx\\
 & \quad \quad + \int_\Omega \int_0^t \left(\frac{1}{p}\,|u(t-s,x)|^p-\frac{1}{p}\,|u(t,x)|^p\right)[-\dot{k}(s)]\,ds\,dx \\
& \quad \quad - \int_\Omega \int_0^t \quad \left(|u(t,x)|^{p-2}u(t,x)[u(t-s,x)-u(t,x)]\right)[-\dot{k}(s)]\,ds\,dx\\
& = \,\frac{1}{p}\,\partial_t(k\ast |u(\cdot)|_p^p)
+\left(-\frac{1}{p}\,|u(t)|_p^p+|u(t)|_p^p\right)k(t)\\
& \quad \quad +\int_0^t \left(\frac{1}{p}\,|u(t-s)|_p^p-\frac{1}{p}\,|u(t)|_p^p+|u(t)|_p^p\right) [-\dot{k}(s)]\,ds\\
& \quad \quad - \int_0^t \left(\int_\Omega \left( |u(t,x)|^{p-2}u(t,x)u(t-s,x)\,dx\right)\right)[-\dot{k}(s)]ds\\
& \ge  \,\partial_t\big(k\ast H(|u(\cdot)|_p)\big) +
\left(-H(|u(t)|_p) + H'(|u(t)|_p)|u(t)|_p\right) k(t)\\
& \quad \quad + \int_0^t \left(H(|u(t-s)|_p)-H(|u(t)|_p)\right)[-\dot{k}(s)]\,ds\\
& \quad \quad - \int_0^t \left(H'(|u(t)|_p)\big(|u(t-s)|_p-|u(t)|_p\big)\right) [-\dot{k}(s)]\,ds\\
& = \,H'(|u(t)|_p)\partial_t\big(k\ast |u(\cdot)|_p)=|u(t)|_p^{p-1}\partial_t\big(k\ast |u(\cdot)|_p\big).
\end{align*}
This proves the lemma. \hfill $\square$

\begin{korollar} \label{lpnorminitial}
Let $1<p<\infty$, $T>0$, and $\Omega$ be an arbitrary measurable subset of $\iR^N$. Let further $k\in
H^1_1([0,T])$ be nonnegative and nonincreasing. Then for any $u_0\in L_p(\Omega)$
and any
$u\in L_p([0,T];L_p(\Omega))$ there holds
\[
|u(t)|_{L_p(\Omega)}^{p-1}\partial_t\left(k\ast\big( |u(\cdot)|_{L_p(\Omega)}-|u_0|_{L_p(\Omega)}\big)\right)(t)\le \int_\Omega |u|^{p-2}u \,\partial_t\big(k\ast [u-u_0]\big)(t)\,dx,
\]
for a.a.\ $t\in (0,T)$.
\end{korollar}
{\em Proof.} By Lemma \ref{lpnorminequ} and H\"older's inequality we have for a.a.\ $t\in (0,T)$
\begin{align*}
\int_\Omega |u|^{p-2}u \,\partial_t\big(k\ast [u-u_0]\big)\,dx  &=
\int_\Omega |u|^{p-2}u \,\partial_t(k\ast u)\,dx-k(t)\int_\Omega |u|^{p-2}u u_0\,dx\,\\
&\ge |u(t)|_p^{p-1}\partial_t\big(k\ast |u(\cdot)|_p\big)(t)
-k(t)|u(t)|_p^{p-1}|u_0|_p\\
&=|u(t)|_p^{p-1}\partial_t\left(k\ast\big( |u(\cdot)|_p-|u_0|_p\big)\right)(t).
\end{align*}
This proves the corollary. \hfill $\square$

${}$

Observe that for sufficiently smooth functions $u$ the statements of Lemma \ref{lpnorminequ} and Corollary \ref{lpnorminitial} remain valid for kernels $k$ which are allowed to be singular at $0$; an important example is given by $k=g_{1-\alpha}$ with $\alpha\in (0,1)$, which leads to an {\em $L_p$-norm inequality for $\partial_t^\alpha$}. Lemma \ref{lpnorminequ} and Corollary \ref{lpnorminitial} also extend to
arbitrary positive measures on $\iR^N$. For illustrative purposes and future reference, let us state a version of Corollary \ref{lpnorminitial}
for functions $u(t)$ taking values in the sequence space $l_p(\iN)$, endowed with the norm $|x|_{l_p}=(\sum_{n=1}^\infty
|x_n|^p)^{1/p}$.
\begin{lemma} \label{folgen}
Let $1<p<\infty$, $T>0$, and $k\in
H^1_1([0,T])$ be nonnegative and nonincreasing. Then for any sequence $u_0=\{(u_0)_n\}_{n=1}^\infty\in l_p(\iN)$
and any function
$u\in L_p([0,T];l_p(\iN))$ there holds with $u(t)=\{(u(t))_n\}_{n=1}^\infty$
\[
|u(t)|_{l_p}^{p-1}\partial_t\left(k\ast\big( |u(\cdot)|_{l_p}-|u_0|_{l_p}\big)\right)(t)\le \sum_{n=1}^\infty |(u(t))_n|^{p-2}(u(t))_n \,\partial_t\big(k\ast [u_n-(u_0)_n]\big)(t),
\]
for a.a.\ $t\in (0,T)$.
\end{lemma}
The next result generalizes the case $p=2$ in Corollary \ref{lpnorminitial} to the Hilbert space setting.
It follows directly from Lemma \ref{fundlemma1}.
\begin{lemma}
Let ${\cal H}$ be a real Hilbert space and $T>0$. Let $k\in
H^1_1([0,T])$ be nonnegative and nonincreasing. Then for any $u_0\in {\cal H}$ and
$u\in L_2([0,T];{\cal H})$ we have
\[
|u(t)|_{{\cal H}}\partial_t\left(k\ast\big( |u(\cdot)|_{{\cal H}}-|u_0|_{{\cal H}}\big)\right)(t)\le
\Big\langle u(t),\frac{d}{dt}\,(k\ast [u-u_0])(t)\Big\rangle_{{\cal H}},\quad \mbox{a.a.}\;t\in (0,T).
\]
\end{lemma}

\section{On the positive part of a subsolution}
\begin{lemma} \label{pospartsubsol}
Let $u_0\in L_2(\Omega)$ and assume that (H) is satisfied. Let
$u\in V(T)$ be a weak subsolution (supersolution) of (\ref{lindiff})--(\ref{lindiff3}) on $\Omega_T$. Then
for any nonnegative $\psi\in \oH^1_2(\Omega)$ there holds
\begin{align*}
\int_\Omega
\Big(\psi\partial_t \big[k_n\ast(u_{+\,(-)}-[u_0]_{+\,(-)})\big]+\big(h_n\ast(AD[u_{+\,(-)}])|D\psi\big)\Big)\,dx\le 0,\quad
\mbox{a.a.}\,t\in (0,T),\,n\in \iN.
\end{align*}
\end{lemma}
{\em Proof.} Suppose $u\in V(T)$ is a weak subsolution of (\ref{lindiff})--(\ref{lindiff3}) on $\Omega_T$. For $\varepsilon>0$, define
\[
H_\varepsilon(y)=\left\{ \begin{array}{l@{\;:\;}l}
(y^2+\varepsilon^2)^{\frac{1}{2}}-\varepsilon & y>0 \\
0 & y\le 0.
\end{array} \right.
\]
Clearly $H_\varepsilon\in C^1(\iR)$ and $H'_\varepsilon\in W^1_\infty(\iR)$. Indeed
\[
H_\varepsilon'(y)=\frac{y}{(y^2+\varepsilon^2)^{\frac{1}{2}}},\quad H''_\varepsilon(y)=\frac{\varepsilon^2}{(y^2+\varepsilon^2)^{\frac{3}{2}}},\quad y>0.
\]
In particular, $H_\varepsilon$ is convex.

Let $\eta\in \oH^{1,1}_2(\Omega_T)\cap L_\infty(\Omega_T)$ be a nonnegative
function with $\eta|_{t=T}=0$. For $t\in (0,T)$ , we take in (\ref{regcond}) the test function $\psi=H'_\varepsilon(u)\eta$, which is admissible,
by boundedness of $H_\varepsilon'$ and $H_\varepsilon''$. We have
\[
D\psi=\eta H''_\varepsilon(u)Du+H'_\varepsilon(u)D\eta,
\]
and thus the resulting inequality can be written as
\begin{align}
\int_\Omega\Big(\eta H'_\varepsilon(u)\partial_t(k_n\ast u)+\big(h_n\ast[ADu]|&\eta H''_\varepsilon(u)Du+H'_\varepsilon(u)D\eta\big)\Big)\,dx\nonumber\\
&\le
\int_\Omega \eta H'_\varepsilon(u)k_n(t)u_0\,dx, \label{esti1}
\end{align}
for a.a.\ $t\in (0,T)$ and any $n\in\iN$.

By the fundamental identity (\ref{fundidentity}) and convexity of $H_\varepsilon$, we have pointwise a.e.\
\begin{align}
H'_\varepsilon(u)\partial_t(k_n\ast u) & \ge \partial_t\big(k_n\ast H_\varepsilon(u)\big)+\big(-H_\varepsilon(u)+H'_\varepsilon(u)u\big)k_n(t)\nonumber\\
& \ge \partial_t\big(k_n\ast H_\varepsilon(u)\big). \label{fund1}
\end{align}
Here we also used convexity of $H_\varepsilon$ to deduce that
\[
-H_\varepsilon(y)+H_\varepsilon'(y)y\ge -H_\varepsilon(0)=0,\quad y\in \iR.
\]

Combining (\ref{esti1}) and (\ref{fund1}), and using that $u_0\le [u_0]_+$, we obtain
\begin{align}
\int_\Omega\Big(\eta\partial_t\big(k_n\ast H_\varepsilon(u)\big) +\big(h_n\ast[ADu]|&\eta H''_\varepsilon(u)Du+H'_\varepsilon(u)D\eta\big)\Big)\,dx\nonumber\\
& \le
\int_\Omega \eta H'_\varepsilon(u)k_n(t)[u_0]_+\,dx, \label{esti2}
\end{align}
for a.a.\ $t\in (0,T)$ and any $n\in\iN$. We integrate (\ref{esti2}) over $(0,T)$,
and then integrate by parts w.r.t.\ time. Sending $n\to \infty$ in the resulting
inequality and using the approximation property of the kernels $h_n$ yields
\begin{align}
\int_0^T\int_\Omega\Big(-\eta_t\big(k\ast H_\varepsilon(u)\big) +\big(ADu|&\eta H''_\varepsilon(u)Du +H'_\varepsilon(u)D\eta\big)\Big)\,dx\,dt \nonumber\\
& \le
\int_0^T\int_\Omega \eta H'_\varepsilon(u)k(t)[u_0]_+\,dx\,dt. \label{esti3}
\end{align}
By positivity of $H''_\varepsilon(u)$ and the parabolicity assumption on $A$, the term $(ADu|\eta H''_\varepsilon(u)Du)$ is nonnegative and thus can be dropped in
(\ref{esti3}). We then send $\varepsilon\to 0$ and use that $H_\varepsilon(y)\to y_+$ and that $H'_\varepsilon(y)\to \chi_{(0,\infty)}(y)$ for all $y\in \iR$, thereby
obtaining
\begin{equation} \label{esti4}
\int_0^T\int_\Omega\big(-\eta_t (k\ast u_+)\big)\,dx\,dt+
\int_0^T\int_{\Omega_+(t)}(ADu|D\eta)\,dx\,dt\le
\int_0^T\int_{\Omega}\eta k(t)[u_0]_+\,dx\,dt,
\end{equation}
where we set $\Omega_+(t)=\{x\in \Omega:\,u(t,x)>0\}$, $t\in (0,T)$. Note that for
any $t\in (0,T)$, we have $D[u(t,x)]_+=Du(t,x)\chi_{\Omega_+(t)}(x)$, a.a.\ $x\in \Omega$. Using this fact and integrating by parts w.r.t.\ time in the integral
on the left-hand side of (\ref{esti4}), it follows that
\begin{equation} \label{esti5}
\int_0^T\int_\Omega\Big(-\eta_t \big(k\ast (u_+-[u_0]_+)\big)+
\big(AD[u_+]|D\eta\big)\Big)\,dx\,dt\le 0,
\end{equation}
for all $\eta\in \oH^{1,1}_2(\Omega_T)\cap L_\infty(\Omega_T)$. By means of an approximation argument that makes use of truncations, it is not difficult to see
that (\ref{esti5}) even holds true for all $\eta\in \oH^{1,1}_2(\Omega_T)$.

Finally, we may argue exactly as in the proof of Lemma 3.1 in \cite{Za} to deduce from (\ref{esti5}) the
assertion in the subsolution case.

If $u$ is a weak supersolution, then $-u$ is a weak subsolution of the same problem with $u_0$ replaced by
$-u_0$, and thus the assertion in the supersolution case follows from the one in the subsolution case.
 \hfill $\square$
\section{Proof of Theorem \ref{result1} and Corollary \ref{korresult1}}
{\em Proof of Theorem \ref{result1}.}
Let $T>0$. Suppose $u\in V(T)$ is a weak subsolution of (\ref{lindiff})--(\ref{lindiff3}) on $\Omega_T$.
Let $\varepsilon>0$ and $\varphi\in C_0^\infty(\Omega)$ be a fixed nonnegative function that is not identically $0$ on $\Omega$. We set
$u_{+,\varepsilon}=u_++\varepsilon \varphi$ and $[u_0]_{+,\varepsilon}=[u_0]_++\varepsilon \varphi$.
Then it follows from Lemma \ref{pospartsubsol}, that for any nonnegative
function $\psi\in \oH^1_2(\Omega)$ and any $n\in\iN$, there holds
\begin{align} \label{esti6}
\int_\Omega
\Big(\psi\partial_t \big(k_n\ast (u_{+,\varepsilon}-[u_0]_{+,\varepsilon})
\big)+\big(h_n\ast[AD[u_+]]|D\psi\big)\Big)\,dx\le
0 ,\quad
\mbox{a.a.}\,t\in (0,T).
\end{align}
For $t\in (0,T)$ we take in (\ref{esti6}) the test function $\psi=u_{+,\varepsilon}$. This gives for any $\varepsilon>0$ and any $n\in \iN$,
\begin{equation} \label{esti7}
\int_\Omega \Big(u_{+,\varepsilon}\partial_t
\big(k_n\ast (u_{+,\varepsilon}-[u_0]_{+,\varepsilon})
\big)+\big(h_n\ast[AD[u_+]]|Du_{+,\varepsilon}\big)\Big)\,dx\le 0 ,\quad
\mbox{a.a.}\,t\in (0,T).
\end{equation}

By Corollary \ref{lpnorminitial} we have for a.a.\
$t\in (0,T)$
\begin{align*}
\int_\Omega u_{+,\varepsilon}\partial_t
\big(k_n\ast (u_{+,\varepsilon}-[u_0]_{+,\varepsilon})
\big)\,dx \ge
\,|u_{+,\varepsilon}(t)|_{L_2(\Omega)}\partial_t \big(k_n\ast (|u_{+,\varepsilon}|_{L_2(\Omega)}-|[u_0]_{+,\varepsilon}|_{L_2(\Omega)})\big).
\end{align*}
Combining (\ref{esti7}) and the previous estimate we arrive at
\begin{align*}
|u_{+,\varepsilon}(t)|_{L_2(\Omega)}\partial_t \big(k_n\ast |u_{+,\varepsilon}|_{L_2(\Omega)}\big)+&
\int_\Omega \big(h_n\ast[AD[u_+]]|Du_{+,\varepsilon}\big)\,dx\\
& \le k_n(t)\big|u_{+,\varepsilon}\big|_{L_2(\Omega)}\,\big|[u_0]_{+,\varepsilon}\big|_{L_2(\Omega)}.
\end{align*}
Since $u_{+,\varepsilon}\ge \varepsilon\varphi\ge 0$ in $\Omega_T$ and $\varphi$ is not identically $0$ in $\Omega$, we have that
\[
W_\epsilon(t):=|u_{+,\varepsilon}(t)|_{L_2(\Omega)}>0,\quad
\mbox{a.a.}\,t\in (0,T),
\]
and thus we obtain for a.a.\ $t\in (0,T)$,
\begin{equation} \label{esti8}
\partial_t \big(k_n\ast W_\epsilon\big)(t)+\,\frac{1}{W_\varepsilon(t)}\,\int_\Omega
\big(h_n\ast[AD[u_+]]|Du_{+,\varepsilon}\big)\,dx
\le k_n(t)\big|[u_0]_{+,\varepsilon}\big|_{L_2(\Omega)},
\end{equation}
for all $\varepsilon>0$ and all $n\in \iN$.

Setting
\[
\mathcal{R}_{\varepsilon,n}(t)=\int_\Omega \big(AD[u_+]-h_n\ast[AD[u_+]]|Du_{+,\varepsilon}\big)\,dx
\]
and using $u_{+,\varepsilon}=u_++\varepsilon \varphi$, (\ref{esti8}) can be rewritten as
\begin{align}
\partial_t \big(k_n\ast W_\epsilon\big)(t)+\,\frac{1}{W_\varepsilon(t)}\,\int_\Omega
\big(AD[u_{+,\varepsilon}]| & Du_{+,\varepsilon}\big)\,dx \le k_n(t)\big|[u_0]_{+,\varepsilon}\big|_{L_2(\Omega)}+\,
\frac{\mathcal{R}_{\varepsilon,n}(t)}{W_\varepsilon(t)}\,\nonumber\\
 & +\,\frac{1}{W_\varepsilon(t)}\int_\Omega
\varepsilon\big(AD\varphi|Du_{+,\varepsilon}\big)\,dx.                                             \label{esti8a}
\end{align}
For any $\delta\in (0,\nu)$ we have
\begin{align*}
\,\frac{1}{W_\varepsilon(t)}\int_\Omega
\varepsilon\big(AD\varphi|Du_{+,\varepsilon}\big)\,dx & \le \,\frac{1}{W_\varepsilon(t)}\,
\left(\frac{\varepsilon^2|AD\varphi|_{L_2(\Omega)}^2}{4\delta}\,+
\,\delta|Du_{+,\varepsilon}|_{L_2(\Omega)}^2\right)\\
& \le \,\frac{\varepsilon|A|_{L_\infty(\Omega_T)}^2
|D\varphi|_{L_2(\Omega)}^2}{4\delta|\varphi|_{L_2(\Omega)}}\,+
\,\frac{\delta|Du_{+,\varepsilon}|_{L_2(\Omega)}^2}{W_\varepsilon(t)}.
\end{align*}
Using this estimate and the parabolicity condition in ($\mathcal{H}$),  we deduce from (\ref{esti8a}) that
\begin{align*}
\partial_t \big(k_n\ast W_\epsilon\big)(t)+(\nu-\delta)\,
\frac{|Du_{+,\varepsilon}(t)|_{L_2(\Omega)}^2}{W_\varepsilon(t)}\le &\, k_n(t)\big|[u_0]_{+,\varepsilon}\big|_{L_2(\Omega)}+\,
\frac{\mathcal{R}_{\varepsilon,n}(t)}{W_\varepsilon(t)}\\
& +\,
\frac{\varepsilon|A|_{L_\infty(\Omega_T)}^2
|D\varphi|_{L_2(\Omega)}^2}{4\delta|\varphi|_{L_2(\Omega)}}.
\end{align*}
By Poincar\'e's inequality, this implies for a.a.\ $t\in (0,T)$ that
\begin{equation} \label{esti8b}
\partial_t \big(k_n\ast W_\epsilon\big)(t)+\lambda_1(\nu-\delta)\,
W_\varepsilon(t)\le \, k_n(t)\big|[u_0]_{+,\varepsilon}\big|_{L_2(\Omega)}+\,
\frac{\mathcal{R}_{\varepsilon,n}(t)}{W_\varepsilon(t)}+\,
\frac{M\varepsilon}{\delta},
\end{equation}
where the positive constant $M=M(|A|_{L_\infty(\Omega_T)},\varphi)$.

Next, denote the right-hand side of (\ref{esti8b}) by $G_{\delta,\varepsilon,n}(t)$
and let $V_{\delta,\varepsilon,n}$ be the solution of the equation
\[
\partial_t \big(k_n\ast V)(t)+\lambda_1(\nu-\delta)\,V(t)= G_{\delta,\varepsilon,n}(t),\quad \mbox{a.a.}\;t\in (0,T),
\]
which exists, since $G_{\delta,\varepsilon,n}\in L_1(0,T)$.
By the comparison principle, see Lemma \ref{compkngl}, we have
\[
W_\varepsilon(t)\le V_{\delta,\varepsilon,n}(t),\quad \mbox{a.a.}\;t\in (0,T),
\]
for all $\varepsilon>0$, $\delta\in (0,\nu)$, and all $n\in \iN$. Sending
$n\to \infty$ and choosing a subsequence if necessary this implies, by Lemma \ref{contdep}, that
\[
W_\varepsilon(t)\le V_{\delta,\varepsilon}(t),\quad \mbox{a.a.}\;t\in (0,T),
\]
for all $\varepsilon>0$, $\delta\in (0,\nu)$, where $V_{\delta,\varepsilon}$ solves
\[
V_{\delta,\varepsilon}+\lambda_1(\nu-\delta)\,l\ast V_{\delta,\varepsilon}=l\ast
G_{\delta,\varepsilon}
\]
and $G_{\delta,\varepsilon,n}\rightarrow G_{\delta,\varepsilon}$ in $L_1([0,T])$, that is
\[
(l\ast G_{\delta,\varepsilon})(t)=\big|[u_0]_{+,\varepsilon}\big|_{L_2(\Omega)}+\,
\frac{M\varepsilon}{\delta}\,(1\ast l)(t),\quad \mbox{a.a.}\;t\in (0,T).
\]
We next send first $\varepsilon\rightarrow 0$ and afterwards $\delta\rightarrow 0$ and choose suitable subsequences.
By the continuous dependence of $V_{\delta,\varepsilon}$ on the parameters $\varepsilon$ and
$\delta$, see e.g.\ \cite{GLS}, and the dominated convergence theorem it follows that
\[
|u_+(t)|_{L_2(\Omega)}\le V(t),\quad \mbox{a.a.}\;t\in (0,T),
\]
where $V$ is the solution of
\[
V+\lambda_1 \nu \,l\ast V=\big|[u_0]_{+}\big|_{L_2(\Omega)}.
\]
Evidently, $V(t)=s_{\nu \lambda_1}(t)|[u_0]_{+}|_{L_2(\Omega)}$, and thus we obtain the desired estimate in the subsolution case, as $T>0$ was arbitrary.

The supersolution case is reduced to the subsolution case by looking at $-u$, which is a global weak subsolution of the same problem with $u_0$ replaced by
$-u_0$. \hfill $\square$

${}$

\noindent{\em Proof of Corollary \ref{korresult1}.} Suppose $u$ is a global weak solution of (\ref{lindiff})--(\ref{lindiff3}). Then from Theorem \ref{result1} we know that the positive and negative part
of $u$, respectively, satisfy the estimate
\begin{equation} \label{kor1est}
\big|u_{+\,(-)}(t,\cdot)\big|_{L_2(\Omega)}\le \big|[u_0]_{+\,(-)}\big|_{L_2(\Omega)}s_{\nu \lambda_1}(t),\quad
\mbox{a.a.}\,t>0.
 \end{equation}
Since $v_+$ and $v_-$ are orthogonal in $L_2(\Omega)$ for any $v\in L_2(\Omega)$, we have
by the Pythagorean theorem that
\[
|u_+(t,\cdot)|^2_{L_2(\Omega)}+|u_-(t,\cdot)|^2_{L_2(\Omega)}=|u(t,\cdot)|^2_{L_2(\Omega)},\;\;
|[u_0]_+|^2_{L_2(\Omega)}+|[u_0]_-|^2_{L_2(\Omega)}=|u_0|^2_{L_2(\Omega)}.
\]
Hence the assertion of Corollary \ref{korresult1} follows from (\ref{kor1est}) by squaring and addition of
the two resulting estimates. \hfill $\square$
\begin{bemerk1}
{\em It is also possible to derive $L_p$-norm decay estimates for suitably defined solutions of (\ref{lindiff})--(\ref{lindiff3}) assuming $u_0\in L_p(\Omega)$, where $1<p<\infty$. In fact, testing the PDE with $|u|^{p-2}u$ and assuming that
$v:=|u|^{(p-2)/2}u\in V(T)$ for all $T>0$ we get with $\rho(p):=4(p-1)/p^2$
\[
|u(t,\cdot)|_p\le s_{\nu \lambda_1 \rho(p)}(t)|u_0|_p,\quad \mbox{a.a.}\;t>0.
\]
Assuming that $u_0$ is bounded and taking the limit $p\to \infty$ gives
$|u(t,\cdot)|_\infty \le |u_0|_\infty$ for a.a.\ $t>0$, since
$\rho(p)\to 0$ and $s_0\equiv 1$. This simple estimate is also a direct consequence of (\ref{maxprin}).
}
\end{bemerk1}

\section{Decay behaviour for some specific examples} \label{kernelexamples}
Theorem \ref{result1} and Corollary \ref{korresult1} show that the decay
properties of the solution to (\ref{lindiff})--(\ref{lindiff3}) is determined
by the behaviour of the relaxation function $s_{\mu}(t)$ with $\mu>0$ for $t\to \infty$.
In this section we discuss in detail this asymptotic behaviour for several examples
of pairs of kernels $(k,l)\in \mathcal{PC}$. We will see that this class of kernels
allows for very different kinds of decay, e.g.\ exponential, algebraic, and logarithmic decay.

We first note that in general $s_{\mu}(t)$ cannot decay faster than the kernel
$k(t)$. Moreover, it is possible that $s_\mu(t)$ does not go to $0$ as $t\to \infty$. In fact we have
\begin{lemma} \label{smuversusk}
Let $(k,l)\in \mathcal{PC}$ and $\mu>0$. Then (i)
\begin{equation} \label{lowerbdd1}
\big[1-s_\mu(t)\big] k(t)\,\le \mu s_\mu(t)\le \,\big[1-s_\mu(t)\big]\,\frac{1}{(1\ast l)(t)} \,,\quad \mbox{a.a.}\;t>0.
\end{equation}
In particular, for any $\delta>0$ there exists $C_\delta>0$ such that
\[
s_\mu(t)\ge \,C_\delta\, k(t),\quad \mbox{a.a.}\;t>\delta.
\]
(ii) $\lim_{t\to \infty} s_\mu(t)=0$ if and only if $l\notin L_{1}(\iR_+)$.
\end{lemma}
{\em Proof.} (i) Recall that $(k,l)\in \mathcal{PC}$ implies that $s_\mu$ and $h_\mu=-\dot{s}_\mu$ are nonnegative, and that
$k_\mu=\mu s_\mu=h_\mu \ast k$. Since $k$ is nonincreasing, it follows that
\begin{align*}
k_\mu(t)\ge (1\ast h_\mu)(t)\, k(t)=(1\ast[-\dot{s}_\mu])(t)\, k(t)=[1-s_\mu(t)\big] k(t),
\end{align*}
which shows the lower bound in (\ref{lowerbdd1}). The upper bound can be easily deduced from the definition of $s_\mu$. In fact, since $s_\mu$ is nonincreasing
we have
\[
s_\mu(t)+\mu s_\mu(t)\,(1\ast l)(t)\le 1,\quad t\ge 0.
\]
The second part of assertion (i) follows directly from (\ref{lowerbdd1}), the monotonicity of $s_\mu$, and the fact that $s_\mu(t)<1$ for any $t>0$.

(ii) By the definition of $s_\mu$ and positivity of $s_\mu$ and $l$, and since
$s_\mu$ is nonincreasing, we have
\[
s_\mu(t)+\mu s_\mu(t) \int_0^t l(\tau)\,d\tau\le 1,\quad t>0,
\]
which implies that $s_\mu(t)\to 0$ as $t\to \infty$ if $l\notin L_{1}(\iR_+)$.
On the other hand, if $l\in L_{1}(\iR_+)$ then
\[
\lim_{t\to \infty} s_\mu(t)=\lim_{z\to 0+} z\hat{s}_\mu(z) =\lim_{z\to 0+} \,
\frac{1}{1+\mu \hat{l}(z)}\,=\,\frac{1}{1+\mu |l|_{L_1(\iR_+)}}\,>0,
\]
by \cite[Theorem 34.3]{Doetsch}. \hfill $\square$
\begin{bemerk1} \label{smubounds}
We point out that (\ref{lowerbdd1}) is equivalent to
\[
\frac{1}{1+\mu \,k(t)^{-1}}\,\le s_\mu(t)\le \,\frac{1}{1+\mu \,(1\ast l)(t)}\,,\quad
\mbox{a.a.}\;t>0.
\]
\end{bemerk1}

${}$

\noindent Before looking at some specific examples we remark that
${\cal PC}$ pairs enjoy a useful stability property with respect
to exponential shifts. Writing $\mathsf{k}_\mu(t)=k(t)e^{-\mu t}$ and $\mathsf{1}_\mu(t)=e^{-\mu t}$, $t>0$,
$\mu\ge 0$ we have
\begin{equation} \label{pcshift}
(k,l)\in {\cal PC} \;\;\Rightarrow\;\; (\mathsf{k}_\mu,\mathsf{l}_\mu+\mu(1\ast
\mathsf{l}_\mu))\in {\cal PC},\;\; \mu\ge 0.
\end{equation}
To see (\ref{pcshift}), observe first that for any $\mu\ge 0$,
$\mathsf{k}_\mu$ is evidently nonnegative and nonincreasing. Multiplying $k\ast l=1$
by $\mathsf{1}_\mu$ gives $\mathsf{k}_\mu \ast \mathsf{l}_\mu=\mathsf{1}_\mu$, which in
turn implies that $\mu \mathsf{k}_\mu\ast 1\ast \mathsf{l}_\mu=\mu 1\ast \mathsf{1}_\mu
=1-\mathsf{1}_\mu$. Adding these relations, we obtain
$\mathsf{k}_\mu\ast[\mathsf{l}_\mu+\mu(1\ast \mathsf{l}_\mu)]=1$.
\begin{bei} \label{beifrac}
The classical time-fractional case. {\em
We consider the pair
\begin{equation} \label{Fbeifrac}
(k,l)=(g_{1-\alpha},g_\alpha),\quad \mbox{where}\;\alpha\in (0,1).
\end{equation}
Recall that the Laplace transform of $g_\beta$, $\beta>0$, is given by $\widehat{g_\beta}(z)=z^{-\beta}$, Re$\,z>0$, and so it is easy to see that $g_{\beta_1}\ast
g_{\beta_2}=g_{\beta_1+\beta_2}$ for all $\beta_1,\beta_2>0$. In particular
$(k,l)\in \mathcal{PC}$.

In the case of (\ref{Fbeifrac})
\[
s_\mu(t)=E_\alpha(-\mu t^\alpha),\quad\mbox{where}\;E_\alpha(z):=\sum_{j=0}^\infty
\,\frac{z^j}{\Gamma(\alpha j+1)}\,,\;z\in \iC,
\]
is the well-known Mittag-Leffler function (see e.g.\ \cite{KST}). Employing the bounds from Remark \ref{smubounds}, a simple computation shows that the Mittag-Leffler function satisfies the estimate
\begin{equation} \label{MLbounds}
\frac{1}{1 + \Gamma (1-\alpha) x} \leq E_{\alpha}(-x) \leq \frac{1}{1+\frac{x}{\Gamma(1+\alpha)}},\,\quad x\ge 0.
\end{equation}
An upper bound of the form $E_\alpha(-x)\le C(\alpha)/(1+x)$, $x\ge 0$, can also be found in \cite{MNV}. Corollary \ref{korresult1} and
(\ref{MLbounds}) yield the {\em algebraic} decay estimate
\[
|u(t,\cdot)|_{L_2(\Omega)}\le \,\frac{1}{1+\Gamma(1-\alpha)^{-1}\nu\lambda_1 t^\alpha}\,|u_0|_{L_2(\Omega)}\le
\,\frac{C(\alpha,\nu,\lambda_1)}{1+t^\alpha}\,|u_0|_{L_2(\Omega)},\quad \mbox{a.a.}\;t>0.
\]
Note that in this example the relaxation function $s_\mu$ has the same decay as the
kernel $k$.
}
\end{bei}
\begin{bei}
The time-fractional case with exponential weight. {\em
We consider
\[
k(t)=g_{1-\alpha}(t)e^{-\gamma t},\quad l(t)=g_{\alpha}(t)e^{-\gamma
t}+\gamma(1\ast[g_{\alpha}e^{-\gamma\cdot}])(t), \quad t>0,
\]
with $\alpha\in (0,1)$ and $\gamma>0$. By the remark prior to Example \ref{beifrac} we have $(k,l)\in \mathcal{PC}$. The Laplace transform of $l$ is given by
\[
\hat{l}(z)=\,\frac{1}{(z+\gamma)^\alpha}\,+
\,\frac{\gamma}{z}\,\frac{1}{(z+\gamma)^\alpha}=\,\frac{(z+\gamma)^{1-\alpha}}{z},
\]
and thus
\[
\widehat{s_\mu}(z)=\,\frac{1}{z}\,\frac{1}{1+\mu \hat{l}(z)}\,=\,\frac{1}{z+\mu(z+\gamma)^{1-\alpha}}.
\]
Let $\omega\in (0,\gamma)$ be the unique solution of $\omega=\mu(\gamma-\omega)^{1-\alpha}$. Note that for fixed $\alpha$ the function $\omega=\omega(\mu,\alpha)$ tends to $0$ as $\mu\to 0$, and $\omega\to\gamma$ as $\mu\to \infty$. Then for any $z\in \iC_+=\{\lambda\in \iC:
\mbox{Re}\,\lambda>0\}$ we have
\[
z-\omega+\mu(z-\omega+\gamma)^{1-\alpha}\neq 0
\]
and thus the Laplace transform of the function $f(t):=s_\mu(t) e^{\omega t}$ is defined for all $z\in \iC_+$. We claim that for some constant $C>0$
\begin{equation} \label{claim1}
|z\hat{f}(z)|\le C \quad \mbox{for all}\; z\in \iC_+.
\end{equation}
Having established (\ref{claim1}) it follows easily that $|z^2 \hat{f}'(z)|$ is
bounded in $\iC_+$ as well. These bounds in turn imply that $f\in L_\infty(\iR_+)$,
by Proposition 0.1 and Corollary 0.1 in \cite{JanI}. Hence $s_\mu(t)\le Me^{-\omega t}$ for all $t\ge 0$. This exponential decay rate is optimal, as $\widehat{s_\mu}$ has a singularity at $-\omega$.

To prove the claim, let $\psi(z):=(z+\gamma-\omega)^{1-\alpha}$ where $|z|<\delta_0:=\gamma-\omega$. There exists $\delta\in (0,\delta_0)$ such that
$\psi(z)=\psi(0)+\psi'(0)z+r(z)$ and the remainder term satisfies
\[
|r(z)|\le \,\frac{1}{2}\,|z|\left(\frac{1}{\mu}+\frac{1-\alpha}{(\gamma-\omega)^{\alpha}}\right)
=\,\frac{1}{2}\,|z|\left(\frac{1}{\mu}+\psi'(0)\right),
\quad \mbox{for all}\,|z|<\delta.
\]
By definition of $\omega$, we have $\mu \psi(0)-\omega=0$, and thus
for $|z|<\delta$, $z\neq 0$, it follows that
\[
|z\hat{f}(z)|= \left|\frac{z}{z+\mu \psi'(0)z+\mu r(z)}\right| = \frac{1}{\mu}\,\left|\frac{1}{\frac{1}{\mu}+\psi'(0)+\frac{r(z)}{z}}\right|
\le \,\frac{2}{1+\mu \psi'(0)}.
\]
On the other hand if $|z|$ is sufficiently large, say $|z|>R$, and Re$\,z\ge 0$ we have
$|z|\ge 2 (\omega + \ \mu |z+\gamma-\omega|^{1-\alpha}$ and thus
\[
|z\hat{f}(z)|\le \frac{|z|}{|z|-\omega-\mu|z+\gamma-\omega|^{1-\alpha}}\le 2.
\]
Finally, by compactness, $z\hat{f}(z)$ is bounded in $\{z\in \overline{\iC_+}:\,\delta\le |z|\le R\}$. This shows (\ref{claim1}).

From Corollary \ref{korresult1} and $s_\mu(t)\le Me^{-\omega t}$ we infer
the {\em exponential} decay estimate
\[
|u(t,\cdot)|_{L_2(\Omega)}\le \,M(\alpha,\nu,\lambda_1)\, e^{-\omega(\nu \lambda_1,\alpha)t}|u_0|_{L_2(\Omega)},\quad \mbox{a.a.}\;t>0.
\]
}
\end{bei}
\begin{bei}
The time-fractional case where $l$ decays exponentially. {\em We consider the situation from the previous example but with the kernels $k$ and $l$ being switched, that is
\[
k(t)=g_{\alpha}(t)e^{-\gamma
t}+\gamma(1\ast[g_{\alpha}e^{-\gamma\cdot}])(t),\quad l(t)=g_{1-\alpha}(t)e^{-\gamma t} \quad t>0.
\]
Note that $\dot{k}(t)=\dot{g}_{\alpha}(t)e^{-\gamma
t}<0$, $t>0$, so that $(k,l)\in \mathcal{PC}$. Since $l\in L_1(\iR_+)$, $s_\mu(t)$
does not go to $0$ as $t\to \infty$, by Lemma \ref{smuversusk}. We have
\[
\lim_{t\to\infty}s_\mu(t)= \,\frac{1}{1+\mu |l|_{L_1(\iR_+)}}=\,\frac{\gamma^{1-\alpha}}{\mu+\gamma^{1-\alpha}}>0.
\]
}
\end{bei}
\begin{bei}
A sum of two fractional derivatives. {\em Let $0<\alpha<\beta<1$ and
\[
k(t)=g_{1-\alpha}(t)+g_{1-\beta}(t), \quad t>0.
\]
Evidently, $k$ is completely monotone, that is, $k$ is in $C^\infty$ and $(-1)^n k^{(n)}(t)\ge 0$ for all $t>0$ and
$n\in \iN\cup \{0\}$. Further $k(0+)=\infty$ and so by Theorem 5.4 in Chapter 5 of \cite{GLS}, the kernel
$k$ has a resolvent $l\in L_{1,loc}(\iR_+)$ of the first kind,
that is $k\ast l=1$ on $(0,\infty)$, and this resolvent is
completely monotone as well. In particular $(k,l)\in {\cal PC}$. The Laplace transforms of $k$ and $l$ are
\[
\hat{k}(z)=\,\frac{1}{z^{1-\alpha}}\,+\,\frac{1}{z^{1-\beta}}\,,\quad \hat{l}(z)=\,\frac{1}{z^\alpha+z^\beta}\,,\quad
z\in \iC_+.
\]
Since $\alpha<\beta$ it is clear that $k(t)\sim g_{1-\alpha}(t)$ as $t\to \infty$. Letting $\mu>0$ we have
\[
\widehat{s_\mu}(z)=\,\frac{1}{z+\mu z \hat{l}(z)}\,=\,\frac{1}{z+\mu\,\frac{z}{z^\alpha+z^\beta}}\,
\sim \,\frac{1}{\mu z^{1-\alpha}}\,\quad\mbox{as}\;z\to 0,
\]
and thus the Karamata-Feller Tauberian theorem, Theorem \ref{Karamata}, implies $\mu s_\mu(t)\sim g_{1-\alpha}(t)$
for $t\to \infty$. We see that the fractional derivative of lower order determines the decay behaviour of the relaxation function $s_\mu$. Observe as well that
\[
\widehat{1\ast l}\,(z)=\,\frac{1}{z}\,\frac{1}{z^\alpha+z^\beta}\,\sim \frac{1}{z^{1+\alpha}}\quad \mbox{as}\;z\to 0,
\]
which yields $(1\ast l)(t)\sim g_{1+\alpha}(t)$ as $t\to \infty$, by Theorem \ref{Karamata}. From this and Remark \ref{smubounds} we infer that there is $T_1>0$ such that
\[
s_\mu(t)\le \,\frac{1}{1+\frac{\mu}{2\Gamma(1+\alpha)}\,t^\alpha},\quad \mbox{for all}\;t\ge T_1, \,\mu>0.
\]
Appealing to Corollary \ref{korresult1} we obtain the decay estimate
\[
|u(t,\cdot)|_{L_2(\Omega)}\le \,\frac{1}{1+\frac{\nu\lambda_1}{2\Gamma(1+\alpha)}\,t^\alpha}\,|u_0|_{L_2(\Omega)}\le
\,\frac{C(\alpha,\nu,\lambda_1)}{1+t^\alpha}\,|u_0|_{L_2(\Omega)},\quad \mbox{a.a.}\;t\ge T_1.
\]
These considerations extend trivally to kernels $k(t)=\sum_{j=1}^m \delta _j g_{1-\alpha_j}(t)$ with $\delta_j>0$ and
$0<\alpha_1<\alpha_2<\ldots<\alpha_m<1$.
}
\end{bei}
\begin{bei}
The distributed order case (ultraslow diffusion). {\em We consider the pair (\ref{ultrapair}) already mentioned in the introduction, that is
\[
k(t)=\int_0^1 g_\beta(t)\,d\beta,\quad l(t)=\int_0^\infty \,\frac{e^{-st}}{1+s}\,{ds},\quad t>0.
\]
Both kernels are nonnegative and nonincreasing. We have
\begin{align*}
\hat{k}(z)=\int_0^1 \widehat{g_\beta}(z)\,d\beta=\int_0^1 z^{-\beta}\,d\beta=
\,\frac{z-1}{z\log z},\quad z\in \iC_+,
\end{align*}
and
\begin{align*}
\hat{l}(z)&=\int_0^\infty e^{-z t} \left(\int_0^\infty \,\frac{e^{-st}}{1+s}\,{ds}\right)\,dt=\int_0^\infty \int_0^\infty e^{-(z+s)t}\,dt \,\frac{1}{1+s}\,ds\\
&=\int_0^\infty \frac{ds}{(z+s)(1+s)}\,=\,\frac{1}{z-1}\,\int_0^\infty
\left(\frac{1}{1+s}-\frac{1}{z+s}\right)\,ds\\
& =\,\frac{1}{z-1}\,\log\left(\frac{1+s}{z+s}\right)\Big|_{s=0}^{s=\infty}\,=\,\frac{\log z}{z-1}\,,\quad z\in \iC_+.
\end{align*}
Thus $(k,l)\in \mathcal{PC}$. The Laplace transform of $s_\mu$ with $\mu>0$ is given by
\[
\widehat{s_\mu}(z)=\,\frac{1}{z}\,\frac{1}{1+\mu \hat{l}(z)}\,=\,\frac{1}{z}\,
\frac{1}{1+\mu \frac{\log z}{z-1}}\,,\quad z\in \iC_+.
\]
We see that $\hat{k}$ and $\mu\widehat{s_\mu}$ have the same asymptotic behaviour near $0$, namely
\[
\hat{k}(z),\,\mu\widehat{s_\mu}(z)\sim \,\frac{1}{z\log(\frac{1}{z})}\quad z\to 0,
\]
where we consider $z>0$. We may apply the Karamata-Feller Tauberian theorem, Theorem \ref{Karamata},
with $L(t):=1/\log t,\,t\ge 2$, $L(t):=1/\log 2,\,t\in (0,2)$, and $\beta=1$, which implies that
\[
k(t),\,\mu s_\mu(t)\sim \,\frac{1}{\log t}\quad t\to \infty.
\]
This can already be found in \cite{Koch08}, where this example is discussed in great detail.

To obtain an estimate that is uniform w.r.t.\ $\mu$ we can also use the upper bound for $s_\mu$ in
Remark \ref{smubounds}. We have $\widehat{1\ast l}(z)=\log(z)/[z(z-1)]$, $z\in \iC_+$,
thus $\widehat{1\ast l}(z)\sim z^{-1}\log(1/z)$ as $z\to 0$, and therefore
$(1\ast l)(t)\sim \log(t)$ as $t\to \infty$, by Theorem \ref{Karamata} with
$L(t)=\log t$ for $t\ge 2$, say, and $\beta=1$. We conclude that there is a number $T_1>1$ (independent of $\mu$) such that $\frac{1}{2}\log t\le (1\ast l)(t)$ for all $t\ge T_1$, and hence
\[
s_\mu(t)\le \,\frac{1}{1+\,\frac{\mu}{2}\,\log t},\quad t\ge T_1.
\]
This together with Corollary \ref{korresult1} yields the {\em logarithmic} decay
estimate
\[
|u(t,\cdot)|_{L_2(\Omega)}\le \,\frac{1}{1+\frac{\nu \lambda_1}{2}\, \log t}\,|u_0|_{L_2(\Omega)},\quad \mbox{a.a.}\;t>T_1.
\]

}
\end{bei}
\begin{bei}
Switching the kernels from the previous example. {\em We consider now the pair
\[
k(t)=\int_0^\infty \,\frac{e^{-st}}{1+s}\,{ds},\quad l(t)=\int_0^1 g_\beta(t)\,d\beta,\quad ,\quad t>0.
\]
From the previous considerations we know already that $(k,l)\in \mathcal{PC}$.
The kernel $k(t)$ in this example behaves like $t^{-1}$ as $t\to \infty$. This can be seen from the representation
\[ k(t)=e^t \int_t^\infty e^{-r}\,\frac{dr}{r},\quad t>0.
\]
In fact, on the one hand we have
\[
k(t)\le \,\frac{e^t}{t}\, \int_t^\infty e^{-r}\,dr=\,\frac{1}{t},
\]
on the other hand we have with $\eta=1+\varepsilon>1$
\[
k(t) \ge e^t \int_t^{\eta t} e^{-r}\,\frac{dr}{r}\ge \,\frac{e^t}{\eta t}\,\left[
-e^{-r}\right]_t^{\eta t}=\,\frac{1}{(1+\varepsilon)t}\,(1-e^{-\varepsilon t}),
\]
and thus $k(t)\ge \frac{1-\varepsilon}{(1+\varepsilon)t}$ for $t>T_\varepsilon$ with sufficiently large $T_\varepsilon$. Since $\varepsilon>0$ is arbitrary, we see that
$k(t)\sim t^{-1}$ as $t\to \infty$.

The Laplace transform of $s_{\mu}$ with $\mu>0$ is given by
\[
\widehat{s_{\mu}}(z) = \,\frac{1}{z}\,\frac{1}{1 + \mu \frac{z-1}{z\log z}}\,=
\,\frac{1}{z+\mu\,\frac{z-1}{\log z}}\,=:\,\frac{1}{\varphi(z)},\quad z\in \iC_+.
\]
Note that $\varphi(z)=z+\mu \int_0^1 z^\beta\,d\beta$, and thus Re$\,\varphi(z)>0$ for all $z\in \overline{\iC_+}\setminus \{0\}$. We see that both $\hat{k}$ and $\mu\widehat{s_\mu}$ behave like $\log(1/z)$ as $z\to 0$. Unfortunately, Theorem \ref{Karamata} does not apply to $s_\mu$ (and $k$) since $\beta=0$ is excluded there. One idea to overcome this
obstacle would be to apply the Karamata-Feller Tauberian theorem to the function $1\ast s_\mu$, which is nondecreasing and has the property that $\widehat{1\ast s_{\mu}}(z)\sim (\mu z)^{-1}\log(1/z)$ as $z\to 0$. This
would show that $(1\ast s_\mu)(t)\sim \mu^{-1} \log t$ as $t\to \infty$. Since $s_\mu$ is
nonincreasing, we would have $t s_{\mu}(t) \leq (1\ast s_{\mu})(t)$ for all $t>0$ and thus
$\mu s_\mu(t)\,\lesssim\, \frac{\log t}{t}$ as $t\to \infty.$ Alternatively, one might look at $1\ast l$. We have
\[
\widehat{1\ast l}\,(z)=\,\frac{z-1}{z^2 \log z}\,\sim \,\frac{1}{z^2 \log(1/z)}\quad \mbox{as}\; z\to 0,
\]
which by Theorem \ref{Karamata} implies that $(1\ast l )(t)\sim t/\log t$ as $t\to \infty$. Remark \ref{smubounds} then
also gives an upper asymptotic estimate for $s_\mu(t)$ as $t\to \infty$ of the form $c\log t/t$.

However, this decay estimate is not optimal. In fact, we will show that $s_\mu$ decays like $c\,t^{-1}$ for any $\mu>0$, that is, the relaxation function has the same algebraic decay as the kernel $k$. To prove the claim, we will show
that the Laplace transform of the function $w(t):=t s_\mu(t)$ satisfies an estimate of the form
\begin{equation} \label{wbounds}
|z\hat{w}(z)|+|z^2\hat{w}'(z)|\le M,\quad \mbox{for all}\;z\in \iC_+
\end{equation}
with some constant $M>0$. Having established (\ref{wbounds}), it follows from Proposition 0.1 and Corollary 0.1 in \cite{JanI} that $w\in L_\infty(\iR_+)$, and thus $s_\mu(t)\le C/t$ for $t>0$.

By a basic property of the Laplace transform we have
\begin{align*}
\hat{w}(z)=\widehat{ts_\mu}\,(z)=-\widehat{s_\mu}'(z)=\frac{\varphi'(z)}{\varphi(z)^2},
\end{align*}
and
\[
\varphi'(z)=1+\mu \,\frac{\log z-1+\frac{1}{z}}{(\log z)^2}\,=1+\mu \int_0^1 \beta z^{\beta-1}\,d\beta.
\]
It is readily seen that as $|z|\to 0$ ($z\in \overline{\iC_+}\setminus \{0\}$) we have $z (\log z)^2 \varphi'(z) \to \mu$ and $\varphi(z)\log z\to -\mu$, and thus $z\hat{w}(z) \to \mu^{-1}$. On the other hand $z\hat{w}(z)\to 0$ as $|z|\to \infty$. By continuity of $z\hat{w}(z)$ in $\overline{\iC_+}\setminus \{0\}$, we thus get an estimate $|z\hat{w}(z)|\le C$ for all $z\in \overline{\iC_+}\setminus \{0\}$. Differentiating once more we obtain
\[
\hat{w}'=\,\frac{\varphi''\varphi-2 (\varphi')^2}{\varphi^3},
\]
with
\[
\varphi''(z)=\mu\,\frac{2z-2-z\log z-\log z}{z^2 (\log z)^3}.
\]
Observe that $z^2(\log z)^2\varphi''(z)\to -\mu$ as $|z|\to 0$. Using this and the above properties and
writing
\[
z^2\hat{w}'(z)=\,\frac{\big[z^2\varphi''(z)(\log z)^2] \big(\varphi(z)\log z\big)-2 \big[\varphi'(z)z(\log z)^2\big]^2\frac{1}{\log z}}{\big(\varphi(z)\log z\big)^3}
\]
we see that $z^2\hat{w}'(z)\to -\mu^{-1}$ as $|z|\to 0$. On the other hand it is not difficult to verify that
$z^2\hat{w}'(z)\to 0$ as $|z|\to \infty$. By continuity of $z^2\hat{w}'(z)$ in $z\in \overline{\iC_+}\setminus \{0\}$,
these observations imply an estimate of the form $|z^2\hat{w}'(z)|\le C_1$ for all $z\in \overline{\iC_+}\setminus \{0\}$. This proves (\ref{wbounds}).

Even more is true. A careful estimation shows that (\ref{wbounds}) holds with some $M$ of the form $M=\frac{\tilde{C}}{\mu}$ where
$\tilde{C}$ is independent of $\mu$. This then leads to an estimate $s_\mu(t)\le \frac{C}{\mu t}$ for all $\mu,\,t>0$.
Since $s_\mu(t)\le 1$ for all $\mu,t\ge 0$, we thus obtain
\begin{equation} \label{smuestexotic}
s_\mu(t)\le \,\frac{c}{1+\mu t},\quad \mbox{for all}\;t,\,\mu\ge 0,
\end{equation}
with some constant $c$ that is independent of $\mu$.

From the previous considerations and Corollary \ref{korresult1} we obtain the algebraic decay estimate
\[
|u(t,\cdot)|_{L_2(\Omega)}\le \,\frac{C(\nu,\lambda_1)}{t}\,|u_0|_{L_2(\Omega)},\quad \mbox{a.a.}\;t>0.
\]
}
\end{bei}

\section{On a basic nonlinear fractional differential equation} \label{nonlin}
Let $\alpha\in (0,1)$, $\gamma,\nu >0$, and $u_0>0$. We are interested in the decay behaviour of the solution to the nonlinear fractional differential equation
\begin{equation} \label{gammagl}
\partial_t^\alpha(u-u_0)+\nu u^\gamma=0,\;\;t\ge 0,\quad\; u(0)=u_0.
\end{equation}

{\em Constructing a subsolution.} Define the positive numbers $\mu$ and $\varepsilon$ by
\[
\mu:= \nu \Gamma(1-\alpha) \Gamma(1+\alpha) u_0^\gamma,
\quad \varepsilon:=\left(\frac{u_0 \Gamma(1+\alpha)}{2\mu}\right)^\frac{1}{\alpha}.
\]
We consider the function
\begin{equation} \label{vdef}
v(t)=\left\{ \begin{array}{l@{\;:\;}l}
u_0-\mu g_{1+\alpha}(t) & t\in [0,\varepsilon] \\
Ct^{-\frac{\alpha}{\gamma}} & t\ge \varepsilon,
\end{array} \right.
\end{equation}
where
\[
C:=\varepsilon^\frac{\alpha}{\gamma} \big(u_0-\mu g_{1+\alpha}(\varepsilon)\big)
=\varepsilon^\frac{\alpha}{\gamma}\,\frac{u_0}{2}.
\]
Observe that $v\in H^1_{1,loc}(\iR_+)$, $v(0)=u_0$, $v(\varepsilon)=u_0/2$, $v$ is nonincreasing, and $v(t)>0$ for all $t\ge 0$.

For $t\in (0,\varepsilon)$ we have
\begin{align*}
\partial_t^\alpha(v-u_0)+\nu v^\gamma &\,= -\mu \partial_t^\alpha g_{1+\alpha}+
\nu (u_0-\mu g_{1+\alpha})^\gamma \\
&\,\le -\mu+\nu u_0^\gamma \le 0,
\end{align*}
by definition of $\mu$. Using $\dot{v}\le 0$ we have for $t> \varepsilon$
\begin{align*}
\partial_t^\alpha(v-u_0)(t) &\,=(g_{1-\alpha}\ast \dot{v})(t)
\le \int_0^\varepsilon g_{1-\alpha}(t-\tau)\dot{v}(\tau)\,d\tau \\
&\, \le g_{1-\alpha}(t)\int_0^\varepsilon \dot{v}(\tau)\,d\tau =-g_{1-\alpha}(t)
\,\frac{u_0}{2}.
\end{align*}
Thus
\begin{align*}
\partial_t^\alpha(v-u_0)+\nu v^\gamma &\,\le -g_{1-\alpha}(t)
\,\frac{u_0}{2}+\nu C^\gamma t^{-\alpha}\\
&\,=-g_{1-\alpha}(t)\Big(\,\frac{u_0}{2}-\nu \Gamma(1-\alpha)\varepsilon^\alpha \left(\frac{u_0}{2}\right)^\gamma\Big)\\
&\, =-g_{1-\alpha}(t)\,\frac{u_0}{2}
\Big(1-\,\frac{\nu}{\mu}\,  \Gamma(1-\alpha) \Gamma(1+\alpha) \left(\frac{u_0}{2}\right)^\gamma\Big)\le 0,
\end{align*}
by definition of $\mu$. Hence $v$ is a subsolution of (\ref{gammagl}).

{\em Constructing a supersolution.} Define $t_0>0$ by means of
\[
t_0^\alpha=\,\frac{u_0^{1-\gamma}}{\nu}\,\left(g_{1-\alpha}\left(\frac{1}{2}\right)+\,
\frac{\alpha}{\gamma}\,\frac{2^{\alpha+\frac{\alpha}{\gamma}}}{\Gamma(2-\alpha)}
\right).
\]
We consider the function
\begin{equation} \label{wdef}
w(t)=\left\{ \begin{array}{l@{\;:\;}l}
u_0 & t\in [0,t_0] \\
Ct^{-\frac{\alpha}{\gamma}} & t\ge t_0,
\end{array} \right. \quad\mbox{with}\;\;C=u_0 t_0^{\frac{\alpha}{\gamma}}.
\end{equation}
For $t<t_0$ we evidently have
\[
\partial_t^\alpha(w-u_0)+\nu w^\gamma=\nu w^\gamma\ge 0.
\]
Next, observe that for $t>t_0$,
\[
\partial_t^\alpha(w-u_0)(t)=(g_{1-\alpha}\ast \dot{w})(t)=-C\,\frac{\alpha}{\gamma}\,\int_{t_0}^t
g_{1-\alpha}(t-\tau)\,\tau^{-\frac{\alpha}{\gamma}-1}\,d\tau.
\]
Assuming $t\in [t_0,2t_0]$ we may thus estimate as follows.
\begin{align*}
\partial_t^\alpha(w-u_0)(t) & \ge -C\,\frac{\alpha}{\gamma}\,t_0^{-\frac{\alpha}{\gamma}-1}g_{2-\alpha}(t-t_0)
\ge -u_0 t_0^{\frac{\alpha}{\gamma}}\,\frac{\alpha}{\gamma}\,t_0^{-\frac{\alpha}{\gamma}-1}\,
\frac{t_0^{1-\alpha}}{\Gamma(2-\alpha)}\\
& \ge -u_0 \,\frac{\alpha}{\gamma \Gamma(2-\alpha)}\,\left(\frac{t}{2}\right)^{-\alpha}
= -\nu w(t)^\gamma\,2^\alpha u_0^{1-\gamma}\,\frac{\alpha}{\nu \gamma \Gamma(2-\alpha) t_0^\alpha}\ge -\nu w(t)^\gamma,
\end{align*}
by definition of $t_0$. For $t>2t_0$ we have
\begin{align*}
\partial_t^\alpha(w-u_0)(t) & =-C\,\frac{\alpha}{\gamma}\, t^{-\alpha-\frac{\alpha}{\gamma}}\int_{t_0/t}^1 g_{1-\alpha}(1-\tau')
\,\tau'^{-\frac{\alpha}{\gamma}-1}\,d\tau'\\
& =-C\,\frac{\alpha}{\gamma}\, t^{-\alpha-\frac{\alpha}{\gamma}}\left(\int_{t_0/t}^{1/2}\ldots+
\int_{1/2}^1\ldots\right)\\
& \ge -u_0 t_0^{\frac{\alpha}{\gamma}}\,\frac{\alpha}{\gamma}\, t^{-\alpha-\frac{\alpha}{\gamma}}
\left(g_{1-\alpha}\left(\frac{1}{2}\right)\,\frac{\gamma}{\alpha}\,
\left(\frac{t_0}{t}\right)^{-\frac{\alpha}{\gamma}}+g_{2-\alpha}\left(\frac{1}{2}\right)
\left(\frac{1}{2}\right)^{-\frac{\alpha}{\gamma}-1} \right)\\
& \ge -\nu w(t)^\gamma \,\frac{u_0^{1-\gamma}}{\nu t_0^\alpha}\,
\left(g_{1-\alpha}\left(\frac{1}{2}\right)+\,\frac{\alpha}{\gamma}\,
\frac{2^{\alpha+\frac{\alpha}{\gamma}}}{\Gamma(2-\alpha)}\right)
= -\nu w(t)^\gamma,
\end{align*}
by the choice of $t_0$. This shows that $w$ is a supersolution of (\ref{gammagl}).

Appealing to Lemma \ref{compstrong} we thus obtain the following result.
\begin{satz} \label{subsuper}
Let $\alpha\in (0,1)$, $\nu, \gamma >0$, and $u_0>0$. Let $u\in H^1_{1,\,loc}(\iR_+)$ be the solution of (\ref{gammagl}) and $v$ and $w$ be defined as in (\ref{vdef}) and
(\ref{wdef}), respectively. Then $v(t)\le u(t)\le w(t)$ for all $t\ge 0$. In particular there exist constants $c_1,\, c_2>0$ such that
\[
\frac{c_1}{1+t^{\frac{\alpha}{\gamma}}}\,\le u(t)\le \, \frac{c_2}{1+t^{\frac{\alpha}{\gamma}}},\quad t\ge 0.
\]
\end{satz}
Theorem \ref{subsuper} shows that the situation in the case $\alpha<1$ differs markedly from that in the case $\alpha=1$,
where we have algebraic decay as $u(t)\thicksim c t^{-1/(\gamma-1)}$ for $\gamma>1$, exponential decay for $\gamma=1$, and extinction in finite time for $\gamma<1$.

\section{On the time-fractional $p$-Laplace equation}
Let $\alpha\in (0,1)$, $1<p<\infty$, and $\Omega\subset \iR^N$ be a bounded Lipschitz domain. We are interested in decay estimates for the solution $u$ of the problem
\begin{align}
\partial_t^\alpha(u-u_0)-\Delta_p u & =0\quad \mbox{in}\;\iR_+\times \Omega,
\nonumber\\
u|_{\partial \Omega} & = 0 \quad \mbox{at}\; \iR_+\times \partial \Omega,
\label{plaplacegl}\\
u|_{t=0} & = u_0 \quad \mbox{in}\; \Omega.\nonumber
\end{align}
Here $\Delta_p u=\mbox{div}\,\big(|Du|^{p-2}Du\big)$. Assuming $u_0\in L_2(\Omega)$
we can define weak solutions of (\ref{plaplacegl}) in a similar way as in the
introduction for (\ref{lindiff})--(\ref{lindiff3}). The natural energy class for
a finite time-interval $[0,T]$ is given by
\begin{align*}
V_p(T):=\{&\,v\in L_{\frac{2}{1-\alpha},\infty}([0,T];L_2(\Omega))\cap
L_p([0,T];\oH^1_p(\Omega))\;
\mbox{such that}\;\\
&\;\;g_{1-\alpha}\ast v\in C([0,T];L_2(\Omega)),
\;\mbox{and}\;(g_{1-\alpha}\ast v)|_{t=0}=0\},
\end{align*}
where the symbol $L_{q,\infty}$ refers to the weak $L_q$-space.
Existence and uniqueness of weak solutions to (\ref{plaplacegl}) in $V_p(T)$ do not seem
to be known in the literature. However we believe that it is possible to
construct weak solutions in $V_p(T)$ using the theory of monotone operators and
the techniques from \cite{ZWH}, at least for $p\ge \frac{2N}{N+2}$. Assuming $u_0\in L_\infty(\Omega)$ global $L_\infty$-bounds for weak solutions have been established
in \cite{VeZa10} by the De Giorgi iteration technique. It is also shown in \cite{VeZa10} that the weak maximum
principle is valid.

In the sequel we write $u\in V_p$ if $u$ belongs to $V_p(T)$ for any $T>0$.
\begin{satz} \label{plaplaceresult}
(i) Suppose that $\frac{2N}{N+2}\le p<\infty$ and that $u_0\in L_2(\Omega)$.
Let $u\in V_p$ be a weak solution of (\ref{plaplacegl}). Then
\begin{equation} \label{plapest1}
|u(t)|_{L_2(\Omega)}\le \frac{C}{1+t^{\frac{\alpha}{p-1}}},\quad \mbox{a.a.}\; t>0,
\end{equation}
where the constant $C=C(\alpha,p,N,\Omega,u_0)$.

(ii) Suppose that $1<p<\frac{2N}{N+2}$, $N> 2$, and that $u_0\in L_s(\Omega)$ where $s=\frac{N(2-p)}{p}$. Let $u\in V_p$ be a weak solution of (\ref{plaplacegl}). Then
\begin{equation}  \label{plapest2}
|u(t)|_{L_s(\Omega)}\le \frac{C}{1+t^{\frac{\alpha}{p-1}}},\quad \mbox{a.a.}\; t>0,
\end{equation}
where the constant $C=C(\alpha,p,N,\Omega,u_0)$.
\end{satz}
{\em Proof.} We proceed by formal a priori estimates. The argument can be made rigorous by adopting the regularization techniques from the proof of Lemma \ref{pospartsubsol} and Theorem \ref{result1}. We may also assume without loss of generality that $u_0$ and $u$ are nonnegative. In fact, by a result analogous to Lemma \ref{pospartsubsol} we may replace $u$ by its positive and negative part, respectively.

(i) In the case $\frac{2N}{N+2}\le p<\infty$ we multiply the PDE by $u$ and integrate
over $\Omega$. This gives
\[
\int_\Omega u \partial_t^\alpha (u-u_0)\,dx +|Du(t)|_p^p\le 0,\quad \mbox{a.a.}\; t\in (0,T).\]
By Corollary \ref{lpnorminitial} and the Sobolev embedding $H^1_p(\Omega)\hookrightarrow L_2(\Omega)$, it follows
that
\[
|u(t)|_2  \partial_t^\alpha \big(|u|_2-|u_0|_2\big)+\nu|u(t)|_2^p\le 0,\quad \mbox{a.a.}\;t\in (0,T),
\]
where $\nu=\nu(\Omega,N,p)$ is a positive constant. Thus $|u(t)|_2$ is a (weak) subsolution of the equation
\[
\partial_t^\alpha(\varphi-\varphi_0)+\nu \varphi^{p-1}=0,\;\,t>0,\quad \varphi(0)=\varphi_0=|u_0|_2.
\]
The desired estimate (\ref{plapest1}) follows now from Lemma \ref{compstrong}, Remark \ref{vergleichweak}, and Theorem \ref{subsuper}.

(ii) We come now to the case $1<p<\frac{2N}{N+2}$. Note first that $s>2$ and that $s\rightarrow 2$ as $p\rightarrow \frac{2N}{N+2}$. Multiplying the PDE by $u^{s-1}$ and integrating
over $\Omega$ yields
\begin{equation} \label{zwischen1}
\int_\Omega u^{s-1} \partial_t^\alpha (u-u_0)\,dx +\mu|Dv(t)|_p^p\le 0,\quad t\in (0,T),
\end{equation}
where
\[
v=u^{\frac{s+p-2}{p}}\quad \mbox{and}\;\;\mu=(s-1)\left(\frac{p}{s+p-2}\right)^p>0.
\]
Corollary \ref{lpnorminitial} implies
\[
\int_\Omega u^{s-1} \partial_t^\alpha (u-u_0)\,dx\ge
|u(t)|_s^{s-1}\partial_t^\alpha \big(|u|_s-|u_0|_s\big).
\]
On the other hand we have the Sobolev embedding $H^1_p(\Omega)\hookrightarrow L_{p^*}(\Omega)$ with
$p^*=\frac{Np}{N-p}$ and the relation $\frac{s+p-2}{p}\cdot p^*=s$, and thus
\[
|u(t)|_s^s=|v(t)|_{p^*}^{p^*}\le C|Dv(t)|_p^{p^*}
\]
for some constant $C>0$. Consequently, it follows from (\ref{zwischen1}) that
\[
|u(t)|_s^{s-1}\partial_t^\alpha \big(|u|_s-|u_0|_s\big)+\,\frac{\mu}{C}\,|u(t)|_s^{\frac{s p}{p^*}}\le 0,
\quad \mbox{a.a.}\;t\in (0,T).
\]
Since $\frac{s p}{p^*}-s+1=p-1$, we may deduce that $|u(t)|_s$ is a (weak) subsolution of
\[
\partial_t^\alpha(\varphi-\varphi_0)+\nu \varphi^{p-1}=0,\;\,t>0,\quad \varphi(0)=\varphi_0=|u_0|_s,\quad \nu=\,\frac{\mu}{C},
\]
 and thus (\ref{plapest2}) follows from Lemma \ref{compstrong}, Remark \ref{vergleichweak}, and Theorem \ref{subsuper}. \hfill $\square$

${}$

\noindent The decay rates in Theorem \ref{plaplaceresult} are optimal, at least for $p> \frac{2N}{N+2}$. In fact, consider a function $u$ of the form
\[
u(t,x)=v(t) w(x),
\]
where $w\in \oH^1_p(\Omega)$ minimizes the functional
\[
F(\psi)=\frac{1}{p}\,\int_\Omega |D\psi|^p\,dx
\]
over the set
\[
K=\{\psi\in \oH^1_p(\Omega):\,|\psi|_{L_2(\Omega)}=1\},
\]
and thus, by standard theory, satisfies $-\Delta_p w=\lambda_1 w$ for some $\lambda_1>0$, cf.\ also \cite[Section 6]{JutLin}. Choosing $v$ to be the solution of
\begin{equation} \label{vvdef}
\partial_t^\alpha (v-1)+\lambda_1 v^{p-1}=0,\;t>0,\quad v(0)=1,
\end{equation}
a short computation then shows that $u$ solves (\ref{plaplacegl}) with $u|_{t=0}=u_0=w$. Theorem \ref{subsuper} implies that $v(t)>0$ for all $t\ge 0$, and thus $u$ cannot become extinct in finite time.
Note that in the case $\alpha=1$ and $p>2$ solutions decay like $c t^{-1/(p-2)}$ (cf.\ \cite{DBUV}) which is not
the same decay rate we obtain when sending $\alpha\rightarrow 1$ in Theorem \ref{plaplaceresult} (ignoring the dependence of $C$ on $\alpha$).
\section{On the time-fractional porous medium equation}
Let $\alpha\in (0,1)$, $0<m<\infty$, and $\Omega\subset \iR^N$ be a bounded Lipschitz domain with $N>2$. This section is devoted to the problem
\begin{align}
\partial_t^\alpha(u-u_0)-\Delta\big(u^m\big) & =0\quad \mbox{in}\;\iR_+\times \Omega,
\nonumber\\
u|_{\partial \Omega} & = 0 \quad \mbox{at}\; \iR_+\times \partial \Omega,
\label{mlaplacegl}\\
u|_{t=0} & = u_0\ge 0 \quad \mbox{in}\; \Omega.\nonumber
\end{align}
We assume that at least $u_0\in L_{m+1}(\Omega)$ and consider nonnegative weak solutions $u$ which are such that for each $T>0$
we have $u^m\in L_2([0,T];\oH^1_2(\Omega))$ and $g_{1-\alpha}\ast (u^{m+1})\in C([0,T];L_2(\Omega))$.
In view of the basic a priori estimates this is a natural class. In the literature nothing seems to be known
on problem (\ref{mlaplacegl}), in paricular existence, uniqueness, and regularity of weak solutions has not been
studied so far.
\begin{satz} \label{mlaplaceresult}
(i) Suppose that $\frac{N-2}{N+2}\le  m<\infty$ and that $u_0\in L_{m+1}(\Omega)$. Let $u$ be a nonnegative weak solution of (\ref{mlaplacegl}).
Then there exists a constant $C=C(\alpha,m,N,\Omega,u_0)$ such that
\[
|u(t)|_{L_{m+1}(\Omega)}\le \frac{C}{1+t^{\frac{\alpha}{m}}},\quad \mbox{a.a.}\; t>0.
\]
(ii) Suppose that $0<m<\frac{N-2}{N+2}$, and that $u_0\in L_s(\Omega)$ where $s=\frac{N(1-m)}{2}$. Let $u$ be a nonnegative weak solution of (\ref{mlaplacegl}). Then there exists a constant $C=C(\alpha,m,N,\Omega,u_0)$ such that
\[
|u(t)|_{L_s(\Omega)}\le \frac{C}{1+t^{\frac{\alpha}{m}}},\quad \mbox{a.a.}\; t>0.
\]
\end{satz}
{\em Proof.} We proceed by formal a priori estimates. The argument can be made rigorous by adopting the regularization techniques from the proof of Theorem \ref{result1}.

(i) Suppose $\frac{N-2}{N+2}\le  m<\infty$. Multiplying the PDE by $u^m$ and integrating over $\Omega$ gives
\[
\int_\Omega u^m \partial_t^\alpha (u-u_0)\,dx +|Dv(t)|_2^2\le 0,\quad \mbox{a.a.}\; t\in (0,T),
\]
where we set $v=u^m$. By H\"older's inequality and Sobolev embedding we have for some constants $c_1,\,c_2>0$ and with $2^*=\frac{2N}{N-2}$
\[
|u(t))|_{m+1}^{2m}\le c_1 |u(t)|_{m\cdot 2^*}^{2m}=c_1 |v(t)|_{2^*}^{2} \le c_2|Dv(t)|_2^2.
\]
Using this and Corollary \ref{lpnorminitial} we obtain the fractional differential inequality
\[
\partial_t^\alpha \big(|u|_{m+1}-|u_0|_{m+1}\big)+\,\frac{1}{c_2}\,|u(t)|_{m+1}^m\le 0,\quad \mbox{a.a.}\; t\in (0,T),
\]
which implies the asserted decay estimate.

(ii) Suppose now that $0<m<\frac{N-2}{N+2}$. Testing the PDE with $u^{s-1}$ we obtain
\[
\int_\Omega u^{s-1} \partial_t^\alpha (u-u_0)\,dx+\mu |Dv(t)|_2^2\le 0,\quad \mbox{a.a.}\; t\in (0,T),
\]
with
\[
v=u^{\frac{s-1+m}{2}}\quad \mbox{and}\quad \mu=\,\frac{4m(s-1)}{(s-1+m)^2}.
\]
Note that $s=2^*\cdot \,\frac{s-1+m}{2}$, and thus by Sobolev embedding we have for some constant $c>0$
\[
|u(t)|_s^{s-1+m}=|v(t)|_{2^*}^2 \le c |Dv(t)|_2^2.
\]
Using this and Corollary \ref{lpnorminitial} we get
\[
\partial_t^\alpha \big(|u|_{s}-|u_0|_{s}\big)+\,\frac{\mu}{c}\,|u(t)|_{s}^m\le 0,\quad \mbox{a.a.}\; t\in (0,T),
\]
which in turn leads to the assertion. \hfill $\square$

${}$

\noindent Recall that in the case $\alpha=1$ and $0<m<1$ any weak solution becomes extinct in finite time provided the initial value $u_0\in L_q(\Omega)$ with $q>1$ and $q\ge s=\frac{N(1-m)}{2}$, cf.\  \cite[Prop.\ 5.23]{VazPM}.
This is no longer the case for the corresponding time-fractional problem, at least for $m>\frac{N-2}{N+2}$. In this case
one can construct similarly as in the previous section a positive solution $u$ of the form $u(t,x)=v(t)w(x)$ where
$v$ solves (\ref{vvdef}) and $w$ is a positive solution of
\[
-\Delta(w^m)=\lambda_1 w\quad \mbox{in}\;\Omega,\quad w=0\;\;\mbox{on}\;\partial\Omega,
\]
with some $\lambda_1>0$; see \cite[Section 4.2]{VazPM} and \cite{BH} for existence of such a $w$. By means of such
a separable solution we also see that the decay rates stated in Theorem \ref{mlaplaceresult} are optimal, at least when $m>\frac{N-2}{N+2}$. Note that in the case $\alpha=1$ and $m>1$ solutions decay like $c t^{-1/(m-1)}$ which is not
the same decay rate we get when sending $\alpha\rightarrow 1$ in Theorem \ref{mlaplaceresult} (ignoring the dependence of $C$ on $\alpha$).


$\mbox{}$

\noindent {\footnotesize {\bf Vicente Vergara}, Universidad de
Tarapac\'{a}, Instituto de Alta Investigaci\'{o}n, Antofagasta N.
1520, Arica, Chile, Email: vvergaraa@uta.cl

$\mbox{}$

\noindent {\bf Rico Zacher}, Martin-Luther-Universit\"at
Halle-Wittenberg, Institut f\"ur Mathematik, Theodor-Lieser-Strasse
5, 06120 Halle, Germany, Email: rico.zacher@mathematik.uni-halle.de

}

\end{document}